\documentclass[a4paper,12pt]{article}
\usepackage{amssymb,amsmath,amsthm}
\usepackage{marvosym}
\usepackage{graphicx}

\usepackage{amsfonts}
\usepackage{latexsym}
\usepackage{color}
\usepackage[english]{babel}
\usepackage{hyperref}
\usepackage{enumitem}
\usepackage{float}
\usepackage{caption}
\usepackage{subcaption}
\numberwithin{equation}{section}

\newtheorem{lemma}{Lemma}[section]
\newtheorem{theorem}[lemma]{Theorem}
\newtheorem{proposition}[lemma]{Proposition}

\newtheorem{remark}[lemma]{Remark}

\theoremstyle{definition}
\newtheorem{manualtheoreminner}{Assumption}
\newenvironment{assumption}[1]{%
	\renewcommand\themanualtheoreminner{#1}%
	\manualtheoreminner
}{\endmanualtheoreminner}
\theoremstyle{plain}

\newtheorem{definition}[lemma]{Definition}
\def\sleq{\lesssim}

\newcommand{\tw}{{\tt w}}
\newcommand{\tf}{{\tt f}}
\newcommand{\tr}{{\tt r}}
\newcommand{\tb}{{\tt b}}

\newcommand{\R}{\mathbb R}
\newcommand{\C}{\mathbb C}

\newcommand{\Z}{\mathbb Z}

\newcommand{\N}{\mathbb N}
\newcommand{\T}{\mathbb T}

\def\norma#1{\|#1  \|}
\def\norm#1{\|#1  \|}

\def\im{{\rm i}}

\def\varep{\varepsilon}

\def\timereg{C^{\infty}_b}

\def\td{{\tt d}}

\newcommand{\csi}{\xi}

\def\cZ{{\mathcal{Z}}}

\def\tk{{\tt k}}

\definecolor{awesome}{rgb}{1.0, 0.13, 0.32}
\definecolor{darkgr}{rgb}{0.0, 0.62, 0.42}
\definecolor{cyan}{rgb}{0.0, 0.72, 0.92}

\def\uno{{\bf 1}}
\def\tM{{\tt M}}

\def\cH{{\mathcal{H}}}
\def\cV{{\mathcal{V}}}

\def\QS{{\mathcal{Q}\kern-0.3pt\mathcal{S}}}
\def\cC{{\mathcal{C}}}
\def\cE{{\mathcal{E}}}

\def\tk{{\tt k}}

\def\tB{{\tt B}}

\def\tR{{\tt R}}

\def\cO{{\mathcal{O}}}
\def\cA{{\mathcal{A}}}
\def\cU{{\mathcal{U}}}

\def\cB{{\mathcal{B}}}
\def\tb{{\tt b}}

\DeclareMathOperator{\OPS}{OPS}

\newcommand{\beq}{\begin{equation}}
\newcommand{\eeq}{\end{equation}}

\begin{document}

\title{{\bf Globally integrable quantum systems and their perturbations}}

\date{}


\author{ Dario Bambusi\footnote{Dipartimento di Matematica, Universit\`a degli Studi di Milano, Via Saldini 50, I-20133
Milano. 
 \textit{Email: } \texttt{dario.bambusi@unimi.it}}, Beatrice Langella\footnote{International School for Advanced Studies (SISSA), via Bonomea 265, I-34136 Trieste.
\textit{Email: } \texttt{beatrice.langella@sissa.it}}.
 }

\maketitle

\begin{abstract}
In this paper we present the notion of globally integrable quantum
system that we introduced in \cite{QN}: we motivate it using the
spectral theory of pseudodifferential operators and then we give some
results on linear and nonlinear perturbations of a globally integrable
quantum system. In particular, we give a spectral result ensuring
stability of most of its eigenvalues under relatively bounded
perturbations, and two results controlling the growth of Sobolev norms
when it is subject either to linear unbounded time dependent
perturbations or a small nonlinear Hamiltonian nonlinear perturbation.
\end{abstract}
\noindent


\medskip



\setcounter{tocdepth}{1}
\tableofcontents
 
\section{Introduction}\label{intro}

In this paper we review some results that we recently got
on perturbations of ``globally integrable quantum
systems''.

One of the key points is the definition of globally integrable quantum
system. Indeed, as far as we know, there is not a generally accepted
definition of integrable quantum system, thus in the first part of the
paper (see Sections \ref{super}, \ref{flows}, \ref{Lie}, \ref{quantumac}) we
first motivate and then give our definition of \emph{globally
integrable
quantum system} (GIQS).  Roughly speaking, a GIQS is a
Hamiltonian operator which is a function of the quantum actions. In
turn we define the quantum actions to be a set of commuting
self-adjoint operators with pure point spectrum contained in
$\Z+\kappa$, where $\kappa\in\R$ is a real number.

The main feature of GIQS is that, on the one hand (as showed by the
result presented in Section \ref{perturbations}) they have very strong
stability properties when perturbed, and on the other hand it turns out
that 
{they cover several nontrivial examples}. That is the reason why we
think they are of interest.

{GQIS exhibit stability properties} both from the spectral point of view
and from a dynamical point of view, when
the Schr\"odinger equation relative to a GIQS is perturbed by a
relatively bounded family of time dependent pseudodifferential
operators, or by a small nonlinear Hamiltonian perturbation.

Such results are reviewed respectively in the
Subsections \ref{linear.pert}
and \ref{nonlinear}. Subsection \ref{linear.pert} is split in a
{paragraph} dealing with the spectral problem and a {paragraph} dealing
with the dynamical problem. In {Paragraph} \ref{spectral} we present
the spectral result, according to which most of the eigenvalues of a
perturbation of a GIQS admit an asymptotic expansion in inverse
powers of their size. In {Paragraph} \ref{growth} we state the
result on time dependent perturbations of the Schr\"odinger equation,
according to which the Sobolev norms of the solution are bounded by 
$|t|^\epsilon$ for all $\epsilon>0$.

Then, in Subsection \ref{nonlinear} we study a nonlinear Hamiltonian
PDE which is a perturbation of the Sch\"odinger equation with
unperturbed Hamiltonian given by a GIQS. Here {we give an almost global existence result}, namely we prove that for initial data
with Sobolev norm of order $\epsilon$ the Sobolev norm of the
corresponding solution remains of order $\epsilon$ for times of order
$\epsilon^{-N}$, $\forall N$. We also give a couple of applications to
concrete nonlinear PDEs. 

The main tool for the proof of all the abstract results is a
clustering property of the spectrum of a GIQS. Such a clustering is
described in Section \ref{clustering} and is constructed by quantizing
the so called  geometric part of the proof of classical Nekhoroshev's
theorem. Such a clustering has also a property that was identified by
Bourgain as a key property for the proof of KAM type results in
PDEs. For this reason we call it Nekhoroshev-Bourgain partition. 

In Section \ref{clustering} we also give a qualitative description of
the way the clustering is used in order to prove the results presented
in Section \ref{perturbations}.

\vskip10pt

We close this introduction by recalling that in the line of research
of this paper one finds a very extensive literature. In particular we
recall that this line of research on the quantization of classical
perturbation theory was initiated many years ago, in particular we
recall \cite{Duclo1, duclos2, combescure, Bel, BG01}, and more
recent developments {were obtained} in \cite{Bam1,Bam2,Bam3} for 1-d systems
and \cite{BGMR1,BGMR2} for some higher dimensional systems. Then we
investigated the case of a nontrivial unperturbed structure of the
resonance relations in \cite{BLMnr,BLMres,BLM_growth,QN}.

The case of nonlinear Hamiltonian {with} semilinear perturbations originated from the
papers \cite{Bam03,BG06} for 1-d systems, the case of the wave
equation on Zoll manifolds was studied in \cite{BDGS} (see
also \cite{BG03,GImekraz}). Recently such a 
theory was extended to the case of quasilinear
perturbations \cite{BD,BMM}. The case of general tori was dealt with
in \cite{BFM22}, which was also the starting point of our
work \cite{QNN}. We also recall the paper \cite{DelIme}, which had a
strong influence on our result \cite{QNN}. 

\vspace{10pt}

\emph{Acknowledgments} The authors would like to thank Francesco
Fass\`{o} for interesting explanations on superintegrable systems and
Giorgio Gubbiotti for a discussion on quantum integrable systems. This project is supported by INDAM and by PRIN 2020 (2020XB3EFL001)  ``Hamiltonian and dispersive PDEs''.

\section{Classical Integrable and Superintegrable Systems. }\label{super}

Let $M$ be an $n$ dimensional Riemannian manifold, and consider the symplectic
manifold given by the cotangent boundle $T^*M$. Let $h\in
C^\infty(T^*M)$ be a Hamiltonian function. Then to $h$ one associates
the corresponding Hamiltonian vector field $X_h$ defined by
$\omega_z(X_h(z),k)=dh(z)k $, $\forall k\in T_z(T^*M)$, where
$\omega_z$ is the symplectic form at $z\in T^*M$. It is well known that
in any canonical coordinate system of $T^*M$ the Hamilton equations $\dot
z=X_h(z)$ take the form
$$
\dot \xi=-\partial_x h\ ,\quad \dot x=\partial_\xi h\ . 
$$
In this case $n=dim(M)$ is called the number of degrees of freedom of the
Hamiltonian system $h$.

One also defines the Poisson Brackets of two functions by
$$
\left\{f_1;f_2\right\}:=d f_1X_{f_2}\ .
$$

\subsection{Integrable systems}

\begin{definition}
  \label{indep}
Let $f_1,...,f_k\in C^\infty(T^*M)$ be $k$ smooth functions. They are
said to be independent if the covectors
$$
df_1(z),...,df_k(z)
$$
are independent for almost every $z$.

The functions are said to be in involution if
$$
\left\{f_j;f_j\right\}=0\ , \forall i,j=1,...,k\ .
$$
\end{definition}

It is well known that the maximal number of functions which are
independent and in involution in a manifold of dimension $n$ is $n$.

\begin{definition}
  \label{integrable}
A Classical Hamiltonian system is said to be integrable if it admits 
$n$ prime integrals which are independent and in involution, namely if
there exist $f_1,...,f_n$ independent and in involution and
s.t. $\left\{h;f_j\right\}\equiv0$. 
\end{definition}

It is well known that the main result of classical integrable systems
is the so called Arnold Liouville theorem, which can be stated as
follows

\begin{theorem}
  \label{arnold}
  Let $h$ be an integrable system; let $z_0\in T^*M$ be s.t. the
  differentials $df_j(z_0)$ are independent and denote
  $m_j=f_j(z_0)$. If the level surface
  $$
S_0:=\left\{z\in T^*M\ :\ f_j(z)=m_j \quad \forall j = 1, \dots, n\right\}
  $$
is compact, then $S_0$ is diffeomorphic to an $n$ dimensional torus and
there exists a neighborhood $\cU$ of $S_0$ foliated in invariant
$n$ dimensional tori. Furthermore there exist $\cA \subseteq \R^n$ and a symplectic coordinate system
$(a,\alpha):\cA\times \T^n\to \cU$, with the property that in such
coordinates $h$ is a function of $a$ only. 
\end{theorem}

{The coordinates $(a,\alpha)$ are called Action Angle coordinates.}

\begin{remark}
  \label{periodiche}
  For the extension to quantum systems, the main property of the
  action coordinates 
is that the coordinate canonically conjugated to
each $a_j$ is an angle. This means that if one considers the
Hamiltonian given by the action $a_j$, such a Hamiltonian generates a
flow which is periodic of period $2\pi$.
\end{remark}

Our main goal is to define the quantum analogue of the action
variable, which in principle is just the operator obtained by
quantizing the classical action. However, quantization is an easy
procedure only for globally defined functions: for this reason we are
particularly interested in situations in which the action variables do
not have singularites, or have only singularities which can be eliminated.

\begin{remark}
  \label{glob}
In this paper, we are interested in the behavior of functions at
infinity. For this reason, the possibly harmful singularities are
only those accumulating at infinity. A typical example are the circular
orbits of the central motion problem, in which the effective radial
Hamiltonian has critical points. It will turn out that such critical
points can be regularized quite easily if they are elliptic, while the
hyperbolic case is not covered by the theory developed here. We also
remark that in more general cases one can meet also singularities of
focus-focus  kind, whose semiclassical theory has been developed
by \cite{PSVG}. 
\end{remark}

\subsection{Superintegrable systems}

A further interesting case where the set of action angle variables is
singular is the case of superintegrable Hamiltonian system, however in
such a case one has that some of the action variables are typically
globally defined and furthermore the Hamiltonian only depends on these
globally defined action variables, a phenomenon which is also related
to Gordon's theorem \cite{Gordon,Nek_Gordon}. As a consequence our theory
will turn out to be applicable also to several superintegrable
systems.

We recall that in the case of superintegrable systems the geometry of
the phase space is described ``semiglobally'' by a classical theorem
by Nekhoroshev \cite{NekInt} (see also \cite{Fasso} and literature
cited therein). We will not state the corresponding theorem,
instead we are going to describe the situation in the case of a free
particle on the 2 dimensional sphere.

So, consider a particle on the two dimensional sphere, which is a
system with two degrees of freedom. Here there exist three independent
integrals of motion which are the three components of the angular
momentum. From these variables one can extract two integrals of motion
which are independent and in involution, namely the $z$ component
$\tM_z$ of the angular momentum and the square $\tM^2$ of the total
angular momentum. A system of action coordinates is given by $\tM$
and by $\tM_z$.  The corresponding
angles can then be constructed by the standard procedure obtaining a
system of action angle variables.

Such a system of coordinates has a singularity when (1) the total
angular momentum is in the direction of the $z$ axis, since here the
two actions are no more independent and (2) when the total angular
momentum vanishes. However, it is clear that the first singularity is
just related to the choice of the $z$ axis and that the other action
$\tM$ is smooth at such a singularity, while at $\tM=0$, there is a
singularity of the action variable, but the Hamiltonian, which is a
function just of $\tM^2$ is smooth also at $\tM=0$. 

This is the typical semiglobal situation described by Nekhoroshev's
theorem \cite{NekInt}: one of the actions is defined on a whole neighborhood
of a level surface of the integrals of motion, while the other one is
not. Typically, a neighborhood $\cU$ of the level surface has the
structure of a bi-fibration $\cU\to^1\cB\to^2 \cA\subset\R$, where
$\cA$ is the domain of the true action on which the Hamiltonian
depends, the fibers of the first map $1$ are one dimensional tori, while
the fibers of the second fibration $2$ are typically compact
surfaces. In the case of the free motion on the sphere, they are locally the
product of a two dimensional sphere and a 1 dimensional
torus. Actually the topological structure of the two dimensional
sphere is the reason why the third action is not globally defined.

\subsection{Conclusion}\label{conc.2}
As a conclusion of this section, we just summarize that classical
integrable systems are systems in which locally the Hamiltonian is a function
of the action variables. In turn, the action variables are functions on
the phase space which, when used as Hamiltonian, generate a flow which
is periodic with period $2\pi$.

Typically the action variables are only defined locally in the phase
space, but there are situations in which they are defined globally.

A particular situation is that of superintegrable systems, in which
the Hamiltonian only depends on a number of actions smaller than the
number of degrees of freedom, and in some interesting situations such
actions are globally defined.

This is a first set of considerations that we use as a basis for the
definition of quantum actions
and quantum integrable system that we will give in Section \ref{def.QN}.

\section{Quantizing Hamiltonians with periodic flow}\label{flows}

As observed in Remark \ref{periodiche} and Subsect. \ref{conc.2}, the
actions of an integrable system generate a periodic Hamiltonian
flow. With this in mind, in this section we discuss quantum operators
obtained by quantizing a Hamiltonian whose flow is periodic. We will
carry on in a parallel way the theory on on $\R^n$ and on compact
manifolds.

To start with consider the case of $\R^n$, where the phase space is
$T^*\R^n\simeq \R^n\times\R^n$. We follow the theory developed by
Helffer-Robert in \cite{HR82}, which is adapted to the study of
anharmonic oscillators in $\R^n$, namely the systems with classical
Hamiltonian
\begin{equation}
  \label{anar}
h_{An}(x, \xi):=\frac{\left\|\xi\right\|^2}{2}+\frac{\left\|x\right\|^{2\ell}}{2\ell}\ ,\quad
\ell\in \N^*\ .
\end{equation}



In order to precisely state the result in \cite{HR82}, we need to introduce some setting.
First we introduce a scale of Hilbert spaces $\{\cH_{An}^s\}_{s \geq
	0}$ adapted to this situation by defining the operator
$K^{(0)}_{An} := (-\Delta + \|x\|^{2\ell})^{\frac {\ell+1}{2\ell}}$
and putting, for $s\geq 0$ $\cH_{An}^s:=D([K_0^{(s)}]^s)$ endowed by the graph norm. For
$s<0$ we also define $\cH_{An}^s$ to be the dual to $\cH_{An}^{-s}$ with respect
to the $L^2$ scalar product.

We also need to introduce the
concept of a quasihomogeneous function. 

\begin{definition}
  \label{quasi}
Fix $\ell\in\N$, then a function $a\in C^{\infty}(\R^{2n}\setminus \{0\})$ is said to be
quasihomogeneous of degree $m \in \R$ if
\begin{equation}
  \label{quasi.1}
a\left(\lambda^{\frac{1}{\ell+1}}x, \lambda^{\frac{\ell}{\ell+1}} \xi\right)=\lambda^{m}a(x, 	\xi) \quad \forall \lambda \in \R^+\,, \quad \forall (x, \xi) \in \R^{2n}\ .
\end{equation}
\end{definition}
Then we define also a class of symbols adapted to this situation:
denote 
            $$\tk_0(x,\xi) := (1+\|x\|^{2\ell}+\|\xi\|^{2})^{\frac{1+\ell}{2 \ell}} \ . $$
\begin{definition}
\label{symbol.ao}
A function $f$ will be called a symbol of order  $m\in\R$ if  $f \in C^\infty(\R^n_x \times \R^n_\xi)$ and 
               $\forall \alpha, \beta \in \N^n$, there exists $C_{\alpha, \beta} >0$ s.t. 
\begin{equation}
\label{es.7}
             \vert \partial_x^\alpha \, \partial_\xi^\beta f( x,\xi)\vert \leq C_{\alpha,\beta} \ \tk_0(x,\xi)^{m-\frac{\ell|\beta| +|\alpha|}{1+\ell}}  \ . 
             \end{equation}
             We will write $f \in S^m_{An}$.
\end{definition}
 To a symbol $f \in S^m_{An}$ we associate the operator
$F\equiv Op^W(f)$ which is obtained by standard Weyl quantization,
 namely
\begin{equation}
	(F u)(x) := \frac{1}{(2\pi)^n} \int_{\R^n} \int_{\R^n} e^{\im \langle x - \eta, \xi \rangle} f\left(\frac{x + \eta}{2}, \xi \right) u(\eta)\, d\eta \ d\xi\,.
\end{equation}
\begin{definition}
  \label{pseudo.an}
  We say that $F$ is a pseudodifferential operator of order $M$ and we
  write $F\in \OPS_{An}^m$ 
   if there
  exist $f \in S^m_{An}$ and a smoothing operator $S$ such that $F =
  Op^W(f) + S$.
  Here by smoothing operator we mean an operator $S$ mapping $\cH_{An}^s$ into $\cH_{An}^{s + N}$ for any $s$ and any $N>0$.
\end{definition}

In this framework we have the following theorem:

\begin{theorem}
  \label{HR82}[Helffer-Robert, \cite{HR82}] Let $a\in S^1_{An}$ be quasihomogeneous
  of degree 1 and assume that all
  the solutions of the corresponding Hamilton equations are periodic
  with period $2\pi$. Assume also that the corresponding Weyl operator
  $A$ is positive definite, then there exist an operator $Q\in
  \OPS_{An}^{-1}$ commuting with $A$ and $\kappa\in\R$ s.t.
  \begin{equation}
    \label{spe}
\sigma(A+Q)\subset \Z+\kappa\ .
  \end{equation}
  By $\sigma(A)$ we are denoting the spectrum of the operator $A$.
\end{theorem}

\begin{figure}
	\centering \includegraphics[scale=0.5]{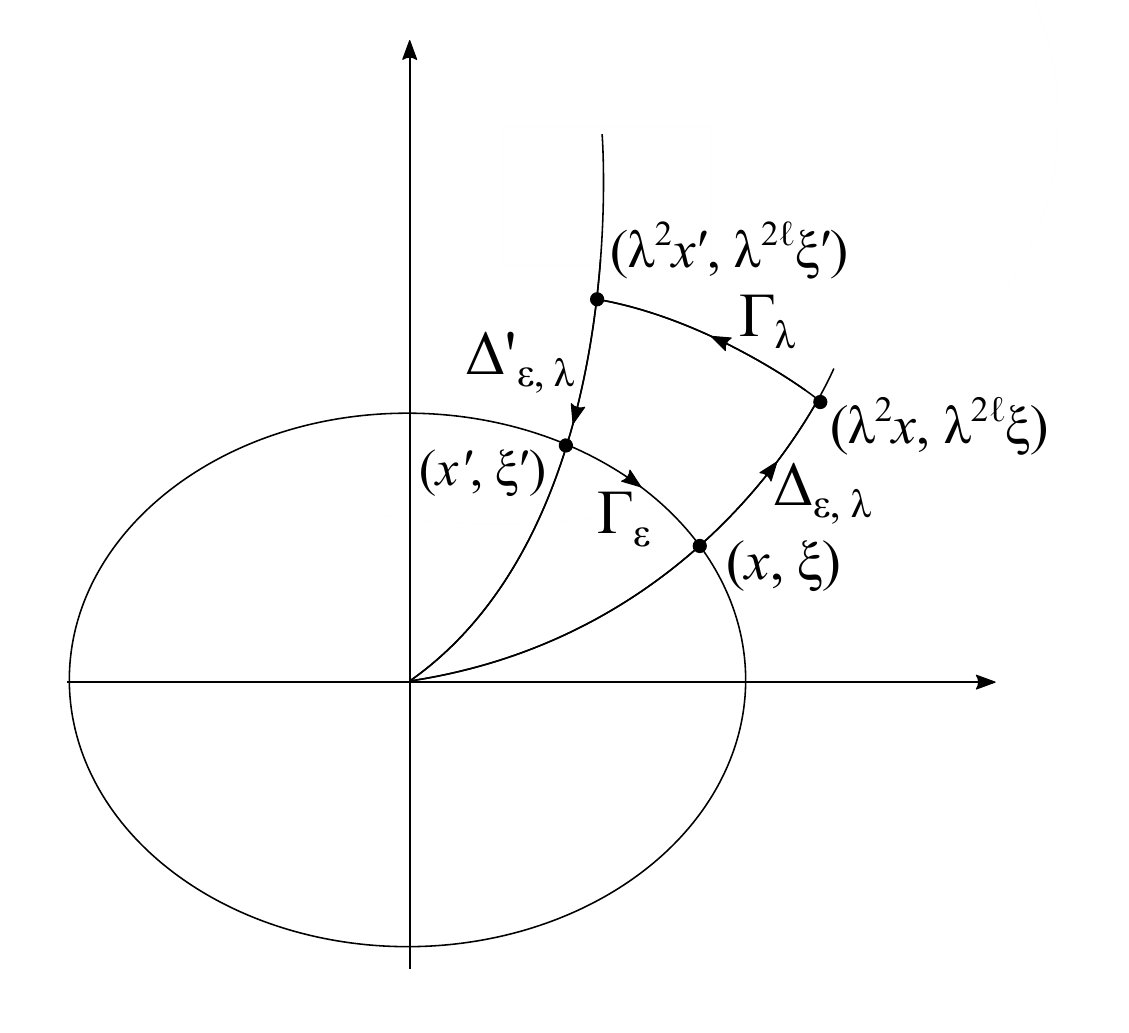}
	\caption{ \small Representation of the path $\gamma$.}
	\label{so_clever}
\end{figure}

\noindent {\it Idea of the proof.} The
idea is first to prove that
\begin{equation}
  \label{quasiuno}
\exp\left(-2\pi \im \left(A-\kappa\right)\right)=\uno+C\,,
\end{equation}
with some $C\in \OPS_{An}^{-1}$ and then to pass to the logarithm. The
intuition is that, by Egorov's theorem, one has that, for any
pseudodifferential operator $B$ one has 
\begin{equation}\label{egorov.teo}
	\exp\left(-2\pi \im A\right)B\exp\left(2\pi \im A\right)=B+l.o.t.
\end{equation}
so that, one expects that, up to lower order terms, $\exp\left(-2\pi \im
A\right)$ is a purely imaginary constant.

The actual proof is slightly different: we are now going to 
summarize how it works. Consider again the operator
$U(2\pi):=\exp(-\im 2\pi A)$. The first step consists in
exploiting Egorov theorem \eqref{egorov.teo} in order to show that,
since the classical flow of $a$ is periodic, $U(2\pi)$ satisfies Beals'
characterization \cite{beals} and therefore it is a pseudodifferential operator in
a suitable class. Then one has to study the principal symbol $u_{2\pi}$ of
$U(2\pi)$. First one gets that it is quasihomogeneous of degree 0. Then
one has to show that it is constant. To this end one defines a
quasihomogeneous function
  $$
M(x,\xi):=\lim_{\lambda\to\infty}u_{2\pi}(\lambda^2x;\lambda^{2\ell}\xi)\ ,
$$
which is easily seen to exist everywhere and finally one proves that
such a function is actually independent of $(x,\xi)$. This last step
is done by writing $u_T(\lambda^2 x, \lambda^{2\ell}\xi)$ as the integral of its differential along the curve
\begin{equation}\label{curva}
	t \mapsto ((\lambda t)^2 x, (\lambda t)^{2\ell}\xi)\,, \quad t\in [0, 1]\,,
\end{equation}
and then by studying the integral of the differential of $u_T$ 
on the 
 closed curve $\gamma$
constructed as in Figure \ref{so_clever}. In particular, given $\varep$ and $\lambda>0$, one takes two points $(x, \xi)$, $(x', \xi')$ on the ellipse $x^{2} +  \xi^{\frac{2}{2 \ell}} = \varep^2$. Then the support of $\gamma$ is built as
$$
\operatorname{\gamma} = \Gamma_\varep \cup \Delta_{\varep, \lambda} \cup \Gamma_{\lambda} \cup \Delta'_{\varep, \lambda}\,,
$$
where
\begin{itemize}
\item[-] $\Gamma_{\varep}$ is the shortest arc of the ellipse joining $(x, \xi)$ and $(x', \xi')$
\item[-] $\Delta_{\varep, \lambda}$ is the arc of the curve \eqref{curva} joining $(x, \xi)$ with $(\lambda^{2} x, \lambda^{2\ell} \xi)$
\item[-] $\Delta'_{\varep, \lambda}$ is constructed as $\Delta_{\varep, \lambda}$ but replacing $(x, \xi)$ with $(x',\xi')$
\item[-] $\Gamma_{\lambda}$ is the arc of ellipse joining $(\lambda^2x, \lambda^{2\ell}\xi)$ with $(\lambda^{2} x', \lambda^{2\ell}\xi')$.
\end{itemize}
Then one proves that, as $\varep \rightarrow 0$ and $\lambda
\rightarrow \infty$, the integrals on $\Gamma_\varep$ and on
$\Gamma_{\lambda}$ tend to $0$. Thus, one gets that the two integrals
on $\Delta_{\varep, \lambda}$ and on $\Delta_{\varep', \lambda'}$
coincide in the limit $\varep \rightarrow 0$ $\lambda \rightarrow
\infty$, namely $M(x, \xi) = M(x', \xi')$. Finally, since $U$ is unitary,
one must have $|M|=1$, and this allows to conclude the proof. \qed

\vskip 10pt

The above theorem clarifies the situation in the case of operators on
$\R^n$. We come now to the case of compact manifolds. So, let $M$ be a
compact Riemannian manifold and denote by $\Psi^m(M)$ the space of
pseudodifferential operators of order $m$ {\it\`{a} la H\"ormander}
(see \cite{hormander}, Chapter 18 for the precise definition). First
we recall that a function $a\in C^{\infty}(T^*M\setminus \{0\})$ is
homogeneous of degree m, if in any canonical coordinate system on
$T^*M$ it has the property that
$$
a(x, \lambda \xi)=\lambda^m a(x, \xi)\ ,\ \forall \lambda>0\ .
$$

The analogue of Theorem \ref{HR82} for the case of compact manifolds
is the following theorem:

\begin{theorem}
  \label{CdV}[Theorem 1.1 of \cite{CdVhelv}]
Let $a\in C^{\infty}(T^*M\setminus \{0\})$ be homogeneous of degree 1
and denote by $A$ its Weyl quantization;
assume that all the orbits of the Hamiltonian vector
field of $a$ are periodic with period $ 2\pi $, then there exists a
pseudodifferential operator $Q\in\Psi^{-1}(M)$ of order $-1$, which
commutes with $A$ and a real number $\kappa$ s.t.
\begin{equation}
  \label{spettro}
\sigma(A+Q)\subset (\Z+\kappa)\ .
\end{equation}
\end{theorem}

Namely, up to a correction smoothing of order $-1$ the quantization of
a Hamiltonian with periodic flow has spectrum contained in a
translation of $\Z$. 

\subsection{Conclusion}
Both in the case of
$\R^n$ and in the case of compact manifolds, the quantization of a
globally defined action is an operator $A$
with the property that there exists an operator $Q$ of order $-1$,
which commutes with $A$ and has spectrum contained in a translation of
$\Z$.

Thus the idea that one can consider is to neglect the regularizing
perturbation and to define directly a quantum
action to be a pseudodifferential operator of order 1 with spectrum
contained in a translation of $\Z$. As shown in the next section, this
has the advantage that in some situations one can avoid to pass through
the classical system.

\section{Two purely quantum examples}\label{Lie}

In this section we give two examples of quantum models where the
Hamiltonian can be described as a function of a few number of operators
$A_1, \dots, A_{\mathtt{d}}$, each of which has spectrum lying on
$\Z+\kappa_j$ for a suitable real $\kappa_j$. 

\vskip 10 pt The first simple example is that of a free quantum particle on
the $d$ dimensional sphere, whose Hamiltonian is the
Laplace-Beltrami operator on the sphere. In this case it is
well known that the spectrum of the Laplacian is given by
$j(j+(n-1))$, and the corresponding eigenspace is spanned by the
spherical Harmonics $Y_{j,l}$. Thus in particular one can write
$-\Delta_g=(A+Q)^2$, with $A$ defined by
$$
AY_{j,l}:=\left(j+\frac{n-1}{2}\right)Y_{j,l}\ ,
$$ and $Q$ defined as a spectral multiplier of order $-1$.

So, in this case one can define a quantum action to be the operator
$A$ and the Laplace-Beltrami operator turns out to be a bounded
perturbation of $A^2$.

The situation is similar on Zoll manifolds, on which one still has
$-\Delta_g=(A+Q)^2$ and $Q$ is a pseudodifferential operator of order
$-1$. 

\vskip 10 pt

A second quantum example is that of a free quantum particle on a
simply connected compact Lie group $G$. We remark that in the case of
$M=SO(3)$ this is a quantum rigid body. Endow $M$ by the bi-invariant
metric, which in the case of the rigid body corresponds to the case of
a spherically symmetric rigid body. The corresponding (quantum)
Hamiltonian is the Laplace-Beltrami operator $-\Delta_g$. As we will
show in a while, in this case the Laplacian turns out to be the
function of a few number of operators with spectrum lying on
$\Z+\kappa$, for some $\kappa \in \R$.

To describe this point we have to rapidly introduce the intrinsic
Fourier series on Lie groups. The key point is that the Fourier
coefficients of a function on a Lie group are labeled by the
irreducible unitary representations of $G$. More precisely, if $\xi$
is an irreducible unitary representation, then the corresponding
Fourier coefficient of a function $\psi$ is an element of the
representation space of such a representation (see for instance
\cite{procesi_padre}). Furthermore one has that such representations
are in 1-1 correspondence with the elements of the cone $\Lambda^+(G)$
of the dominant weights, defined by
\begin{equation}\label{cono.pesi}
\Lambda^+(G) = \left\{ \tw \in \R^d\ \left|\ \tw = \tw^1 \tf_1 + \dots
+ \tw^d \tf_d\,, \quad \tw^j \in \N\ , \  \forall j= 1, \dots, d \right.\right\}\,,
\end{equation}
where $\tf_1, \dots, \tf_d \in \R^d$ are the \emph{fundamental weights} of
$G$. In this language the Laplacian is a Fourier
multiplier. Precisely, given a dominant weight, $\tw$ consider its
decomposition on the basis $\tf_j$, namely $\tw=\sum_{j=1}^d\tw^j
\tf_j$, then the Laplacian is the Fourier multiplier by
$$
\sum_{j,i}(\tw^j+1)(\tw^i+1)\tf_i\cdot\tf_j-\sum_{j,i}\tf_i\cdot\tf_j\ ,
$$
which is a homogeneous quadratic function of the basic Fourier
multipliers
\begin{equation}
  \label{az.Lie}
A_j:=\text{multiplier\ by}\ (\tw^j+1)\ .
  \end{equation}
\begin{remark}
  \label{Lie.2}
Since $\tw^j$ are integers, the spectrum of such operators are just
contained in $\Z$. Furthermore, exploiting the results of
\cite{Fischer}, in \cite{QN} we proved that such operators are
actually pseudodifferential operators in the sense of H\"ormander.
\end{remark}

Actually it can be seen that the situation is absolutely similar in homogeneous spaces. We also remark that the sphere is also a
homogeneous space, and this suggests that there should be a connection
between the two examples of this section. However, we also remark that
Zoll manifolds are not, in general, homogeneous spaces.

\section{Globally integrable quantum systems}\label{quantumac}

The theorems and the examples of the previous sections give a strong
hint on a possible definition of quantum action and quantum integrable
system.

We will treat in a unified way both the situation of a compact manifold
or that of $\R^n$. In particular on $\R^n$ we are interested to the
anharmonic oscillator with quantum Hamiltonian
\begin{equation}
  \label{ana}
H_{An}:=-\frac{\Delta}{2}+\frac{\left\|x\right\|^{2\ell}}{2\ell}\ ,
\end{equation}
with some fixed $\ell\in\N$, $\ell\geq 1$.  We define
\begin{equation}
  \label{k0}
  K_0:=\left\{
  \begin{matrix}
    (\uno-\Delta_g)^{1/2}\ &\text{if }& M & \text{is\ compact}
    \\
    (\uno+H_{An})^{(\ell+1)/2\ell}\ &\text{if }& M & \text{is}\ \R^{n}  
  \end{matrix}
  \right.\,.
  \end{equation}
Then we consider the scale of Hilbert
spaces $\cH^s:=D(K_0^s)$, if $s>0$ and $\cH^s$ the dual of $\cH^{-s}$
with respect to the $L^2$ scalar product, if $s<0$.

Concerning pseudodifferential operators, we define
\begin{equation}
  \label{psido}
  \cA^m:=\left\{
  \begin{matrix}
    \Psi^m(M)\ &\text{if }& M & \text{is\ compact}
    \\
\OPS_{An}^m\ &\text{if }& M & \text{is}\ \R^{n}  
  \end{matrix}
  \right.
  \end{equation}
  
\subsection{Definitions}\label{def.QN}
As anticipated above, the actions will be a set of
commuting self-adjoint operators, each of which with spectrum contained in a
translation of $\Z$. 

Before giving the precise definition we recall that, given a set of
self-adjoint operators $A_1,...,A_d$ with pure point spectrum, their
joint spectrum is defined as follows:
\begin{definition}
  \label{joint}
 {The joint spectrum $\Lambda$ of the operators $A_j$ is
the set of the ${a = (a^1, \dots, a^d)}\in\R^d$ s.t. there
exists $\psi\in\cH$ with 
	\begin{equation}
	\label{lambda}
        A_j\psi=a^j\psi
       \ ,\quad \forall j=1,...,d\ .
	\end{equation}
}
\end{definition}

\begin{definition}[Set of quantum actions]  \label{actions} A set of self-adjoint
  pseudodifferential operators $(A_1,...,A_d)$, with $A_j\in\cA^1$ are
  said to be a set of quantum actions if  
  \begin{itemize}
  \item[i.] They are pairwise commuting
  \item[ii.] $\exists c_1>0$ s.t. $c_1 K_0^2<
		\uno+\sum_{j=0}^d A_j^2$.
		\item[iii.] There exist a convex closed cone $\cC\subseteq \R^{d }$
		and a vector
		$
		\kappa=(\kappa_1,...,\kappa_{d })\in\R^d
		$,
		such that the joint spectrum $\Lambda$ of the $A_j$'s
		fulfills
		\begin{equation}
		\label{I1}
		\Lambda\subset (\Z^{d}+\kappa)\cap \cC\ .
		\end{equation}
  \end{itemize}
\end{definition}

\begin{definition}
  \label{qi}
A self-adjoint operator $H_L$ is said to be \emph{globally integrable} if
there exists a function $h_L:\R^d\to\R$ with
the property that
\begin{equation}
  \label{inte}
H_L=h_L(A_1,...,A_d)\ ,
\end{equation}
with the operator function spectrally defined. 
\end{definition}

We remark that Definition \ref{qi} is given in order to cover also the
case of the quantization of superintegrable systems, in which the
total number of globally defined action variables is less than the
number of degrees of freedom. In particular with this definition one
renounces to any assumption ensuring ``completeness'' of the set of
quantum actions.
This is recovered in the case where
the multiplicity of the common eigenvalues of the actions is 1. Indeed
one has the following very simple remark

\begin{remark}
  \label{commu}
Assume that the multiplicity of the common eigenfunctions of the
action is 1, namely that $\forall a\in\Lambda $ there exists a {\bf unique}
$\psi_a\in \cH$ s.t. $\left\|\psi\right\|=1$ and \eqref{lambda} holds,
then if $F\in\cA^m$ for some $m$  is a self-adjoint
operator on $\cH$ with the property that
$$
[A_j,F]=0\ ,\quad \forall j=1,...,d
$$ then there exists a function $f:\Lambda\to\R$ s.t. $F=f(A)$. Indeed one
can just define $f(a):=\langle\psi_a;F\psi_a\rangle_{\cH}$.
\end{remark}

The Definitions \ref{actions} and \ref{qi} of system of quantum actions and of integrable quantum
system we gave have a strong limitation. This is hidden in the statements
of Theorems \ref{HR82} and \ref{CdV}. Indeed such theorems give the
spectrum of an operator obtained by quantizing a function $a$ which is
\emph{globally defined} on the phase space. This is the reason why we
added the word \emph{globally} to the word integrable in the above
definition. Furthermore, in view of this limitation it is important to
present some nontrivial examples of \emph{globally integrable quantum
  systems}.

\subsection{Examples of globally integrable quantum systems}\label{EGIQS}

As proven in \cite{QN}, some examples of globally integrable quantum systems are:
\begin{enumerate}
\item $(M, g) = (\T^d, g)$, with $g$ an arbitrary flat metric, and
  $H_L = -\Delta_g$. Here the quantum actions are given by $A_j := \im
  \partial_j$, and one easily sees that $-\Delta_g = H_L =
  \sum_{i,j=1}^d g^{ij}A_iA_j$, where $g^{ij}$ is the inverse of the matrix of
  the metric.

\item $(M,g)$ a Zoll manifold. This is the situation studied in
  \cite{CdVhelv}, and is the same one described for spheres in Section
  \ref{Lie}. As for spheres, one has $-\Delta_g = (A + Q)^2$, for some
  commuting operators $Q \in \Psi^{-1}$ and $A \in \Psi^1$. The
  quantum action is $A$, and one again has $H_L = A^2 =-\Delta_g +
  \text{lower order terms}$.
  
\item $(M, g)$ a rotation invariant surface as follows: $M$ is the
  level surface $f(x,y,z) = 1$ of a $C^\infty$ function $f:\R^3
  \rightarrow \R$ which is a submersion at $f(x,y,z) = 1$ and invariant under rotation around the $z$ axis.
We  endowed $M$
  with metric $g$ induced by the embedding. If $\theta \in [0, L]$ is the
  curvilinear abscissa along the geodesic $\gamma$ as in Figure
  \ref{fig.bottiglia}, and we introduce cylindrical coordinates $(r,
  \theta, \phi)$, then in such coordinates one has $M= \{(r(\theta),
  \theta, \phi)\ |\ \theta\in [0, L]\,, \quad \phi \in \T\}$. We
  assume that the function $r(\theta)$ has a unique critical point in
  $(0, L)$. This case was covered by \cite{CdV}, where it was shown
  that there exists a GIQS $H_L$ which coincides with  $-\Delta_g$ up
  to lower order corrections. Here the quantum actions are a
  lower order perturbation of of the quantization of the actions of
  the principal symbol of $H_L$.
\begin{figure}
		\centering \includegraphics[scale=0.4]{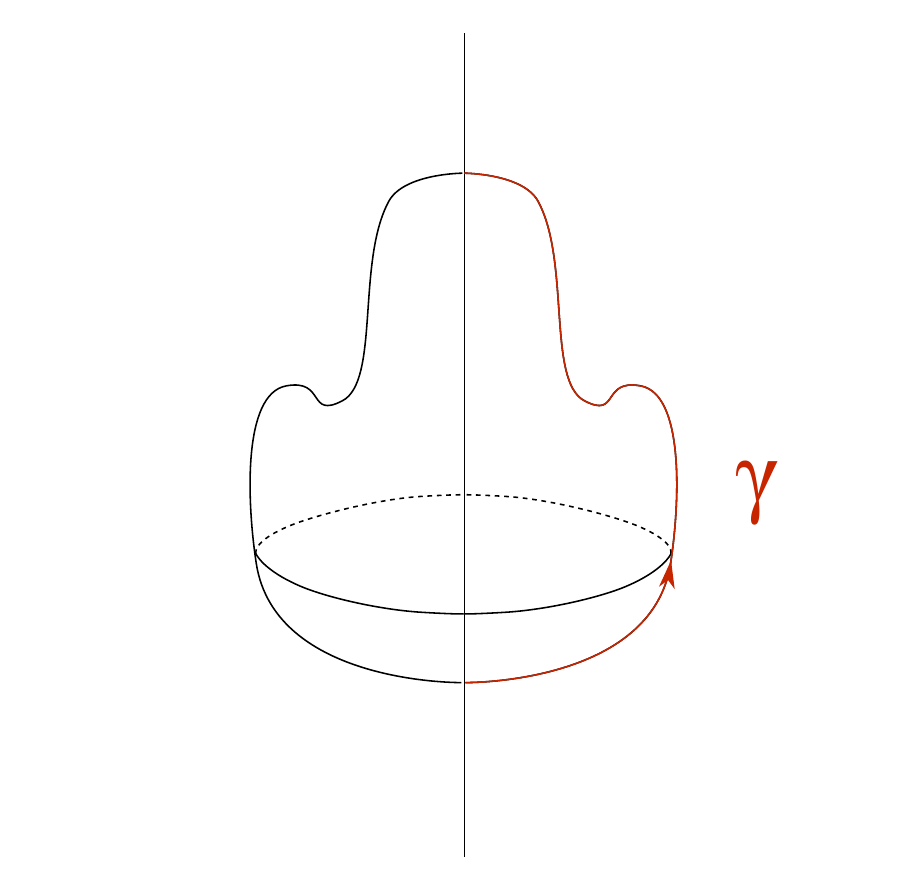}
	\caption{ \small An admissible rotation invariant surface $M$. In red, the geodesic $\gamma$ parametrized by $\theta \in [0, L]$.}
	\label{fig.bottiglia}
\end{figure}
\item $(M, g)$ a simply connected, compact Lie Group, with $H_L$ a bounded perturbation of $-\Delta_g$. This is the situation described in Section \ref{Lie}; further details will be discussed in Section \ref{sez.per.esempio}.

\item $(M, g) = (\R^2, g)$ with $g$ the Euclidean metric and $H_L$ a
  correction of $H_{An}$, cf. Eq. \eqref{ana}. This model will be discussed in detail in Section \ref{sez.per.esempio}
\end{enumerate}
\begin{remark}
  We point out that, in the examples 1,3 and 5, any $a$ in the joint
  spectrum $\Lambda$ has multiplicity 1 as a joint eigenvalue for
  $A_1, \dots, A_d$, so we are in the situation of Remark \ref{commu}.

  On the contrary, in Examples 2 and 4 one deals with superintegrable
  quantum systems.
\end{remark}
\subsection{Semiclassical integrable systems}
We end this section by recalling the definition of a semiclassical
integrable system given by Colin de Verdi\`{e}re in his book
\cite{CdVBook}. First we recall that Colin de Verdi\`{e}re deals with
\emph{semiclassical Hamiltonian}\footnote{for a precise definition, we
  refer to \cite{CdVBook}} operators on compact manifolds and studies
the behavior of the various objects as $\hbar\to0$.

\begin{definition}
  \label{CdV.def}
A semiclassical Hamiltonian $H$ on an $n$ dimensional manifold is said
to be \emph{semiclassically completely integrable} if there exist $n$
selfadjoint semiclassical operators $F_1,...,F_n$ of order $0$ which
are pairwise commuting s.t.
\begin{itemize}
\item[i.] The principal symbols $f_1,...,f_n$ are in involution and
  almost everywhere independent;
\item[ii.] There exists a function $h_L\in C^\infty(\R^n:\R)$ s.t.
  $$
H=h_L(F_1,...,F_n)\ .
  $$
\end{itemize}
\end{definition}

We remark that this definition does not make any reference to the
action angle variables, but just to the existence of a sufficient number
of ``prime integrals'' of the system. 

In \cite{CdVBook} one can also find a semiclassical analogue of Arnold
Liouville theorem, showing that microlocally close to an invariant
nondegenerate torus of the classical system, one can conjugate the
Hamiltonian to a function of the operators $-i\hbar\partial_{j}$ (see
also Charbonnel \cite{ChaIsa}) on $\T^d$. We emphasize that the
semiclassically integrable systems dealt with in such studies are more
general than those covered by our definition, but on the contrary
these semiclassical studies only give local results (see also
\cite{PSVG} for some semiglobal results). We remark that we do not
know if the results of Sect. \ref{perturbations} can be extended to
semiclassical integrable systems.

\section{Perturbations of globally integrable quantum
  systems}\label{perturbations} 

In this section we will give the statement of two kind of results on
perturbations of globally integrable quantum systems, the first one
deals with linear perturbations, while the second deals with small nonlinear
Hamiltonian perturbations. 

The main assumption on the GIQS that we perturb is that $h_L$ is
homogeneous and steep. The main tool used in the proof is a result of
clustering of the joint spectrum $\Lambda$ that is
given in the forthcoming Section \ref{clustering} (see Theorem
\ref{teo.partizione}). 

\vskip 5pt

We recall that homogeneous functions are typically singular at the
origin, but, in the context of pseudodifferential operators, the
behavior of functions in a neighborhood of the origin is not
important. In particular, in order to get rid of singularities at the
origin remaining in the class of symbols, we give the following:
\begin{definition}
	A function $f \in C^\infty(\R^d)$ is said to be homogeneous of
        degree $m$ at infinity if there exists $R>0$ such that
	$$
	f(\lambda a) = \lambda^{m} f(a) \quad \forall \lambda >1\,, \quad \forall |a| \geq R\,.
	$$
\end{definition}

     We recall, from \cite{GCB}, the definition of steepness:
\begin{definition}[Steepness]	\label{steep.def} 
	Let $\cU\subset \R^d$ be a bounded connected open set with nonempty
	interior.  A function $h_0 \in C^ 1 (\cU )$, is said to be steep in
	$\cU$ with {steepness radius $\tr$,} steepness indices $\alpha_1,\dots
	,\alpha_{d-1}$ and (strictly positive) steepness coefficients $\tB_ 1
	, . . . , \tB_{d-1}$, if its gradient $w_i(a):=\frac{\partial
		h_0}{\partial a_i}(a)$ fulfills:
	$\displaystyle{\inf_{a\in\cU}\norma{w(a)}>0}$ and for any
	$a\in\cU$ and for any $s$ dimensional linear subspace $M\subset \R^d$
	orthogonal to $w(a)$, one has
	\begin{equation}
	\label{steep}
	\max_{0\leq\eta\leq\xi}\min_{u\in M:\norm u=1}\norma{\Pi_M
		w(a+\eta u)}\geq \tB_s\xi^{\alpha_s} \quad \forall \xi \in (0, \tr]\ ,
	\end{equation}
	where $\Pi_M$ is the orthogonal projector on $M$; the quantities $u$
	and $\eta$ are also subject to the limitation $a+\eta u\in\cU$. 
\end{definition}

It is well known that steepness is generic. We also recall that all
convex or quasiconvex functions are also steep.

A very useful  sufficient condition for steepness is contained in the
following theorem by Niederman \cite{Nied06}. This is the theorem that
we used in order to verify steepness in our concrete applications.

\begin{theorem}[Niederman] \label{nie06} Let $h_L$ be a function real analytic
        in an open set $\cU\subset\R^d$. Assume that {$h_L$ has no
          critical points in $\cU$,} and its restriction $h\big|_{M}$
        to any affine subspace {$M\subset \R^d$} admits only isolated
        critical points, then $h_L$ is steep on any
        compact subset of $\cU$.
\end{theorem}

So, we concentrate on perturbations of linear systems with Hamiltonian
$H_L$ which is globally integrable quantum system, namely $H_L=h_L(A)$
and make the following assumption

\begin{assumption}{(L)}\label{L}
$h_L$ is homogeneous of degree $\td>1$ at infinity and steep
        in an open set $\cU\subset B_{1}\setminus B_{1/2}$ with the
        property that
        $$
\Lambda\cap \cU=\Lambda\cap (B_{1}\setminus B_{1/2})\ .
        $$
\end{assumption}

\subsection{Linear Perturbations: spectral stability and growth of
  Sobolev norms} \label{linear.pert}

\subsubsection{Stability of the majority of eigenvalues}\label{spectral}

The first result we are going to present is a stability result of the
majority of eigenvalues of $H_L$ when one adds a perturbation which is
relatively bounded with respect to $H_L$ this is a generalization of
Theorems 2.8 of \cite{BLMnr} and 2.8 of \cite{anarmonico}.

First we consider the general case where the common eigenvalue possible
have some multiplicity.

In this case it is possible to compute the first correction to the
eigenvalues bifurcating from the majority of eigenvalues of
$H_L$. This requires to define the average of the perturbation $V$,
namely the operator
\begin{equation}
  \label{average}
\langle V\rangle:=\frac{1}{(2\pi)^d}\int_{\T^d}e^{-\im\varphi\cdot
  A}Ve^{\im\varphi\cdot A}d \varphi\,,
\end{equation}
which commutes with all the actions. Denote also
\begin{equation}\label{Pi.a}
	\Pi_{a}:=\Pi_{a^1}...\Pi_{a^d}
\end{equation}
the projector on the common
eigenspace of the actions corresponding to the common eigenvalue
$a=(a^1,...,a^d)$, and denote by $n_a$ the dimension of $\Pi_a\cH$,
namely the multiplicity of the common eigenvalues, then 
the operator $\Pi_a\langle
V\rangle\Pi_a$, gives the first variation of the eigenvalues
bifurcating from $\lambda_a$. Precisely, denote by 
$\mu_{a,j}$ the
eigenvalues of $\Pi_a\langle
V\rangle\Pi_a$, then we have the following:
\begin{theorem}
  \label{rilinea}
  	Consider the operator
	\begin{equation}
	\label{ope1}
	H:=H_L+V \ ,
	\end{equation}
	assume that $H_L$ fulfills assumption \ref{L} and
        $V\in\cA^{\tb}$ with $\tb<\td$. Then there exist a set $\Omega
        \subset \R^d$ and number $m_2<\tb$ such that the following
        holds. $\Omega \cap \Lambda$ has density one at infinity,
        namely there exists $\rho>0$ such that
	\begin{equation}
	\label{density}
	1-\frac{\sharp(\Omega\cap\Lambda \cap B_R)}{\sharp(B_R \cap
		\Lambda)} = \cO(R^{-\rho})\, ,
	\end{equation}
and $\forall a\in\Omega$ and any
$j=1,...,n_a$ there exists an eigenvalue $\lambda_{a,j}$ of
\eqref{ope1} fulfilling
$$
|\lambda_{a,j}-\lambda_a-\mu_{a,j}|\sleq \langle a\rangle^{m_2}\ . 
$$
\end{theorem}

In the case in which the common eigenvalues of the actions have no
multiplicity one can give a more precise result, to this end we need
first the following:
\begin{definition}
	Given a function $f: \cV \subseteq \R^d \rightarrow \R$, we say that $f$ admits an asymptotic expansion,
	\begin{equation}
		f \sim \sum_{j} f_j\,,
	\end{equation}
	if there exist a real sequence $\{m_j\}_{j \in \N}$ and a sequence $\{f_j\}_{j \in \N}$ of smooth functions $f_j: \cV \rightarrow \R$ such that
	\begin{gather}
		m_{j+1} < m_{j} \quad \forall j\,, \quad m_j \rightarrow -\infty \quad \text{as} \quad j \rightarrow \infty\,,\\
		|f_j(a)| \leq C_j \langle a\rangle^{m_j} \quad \forall a\,,\\
		\left|f(a) - \sum_{j=0}^{N} f_j(a)\right|\leq C_N \langle a \rangle^{m_{N+1}} \quad \forall a, \quad \forall N\in \N\,.
	\end{gather}
\end{definition}

\begin{theorem}
	\label{maint}
Under the same assumptions of Theorem \ref{rilinea}, assume also that the multiplicity of each common
	eigenvalues of the actions is 1 (as in Remark \ref{commu}). Then
	there exist a sequence of smooth functions $z_j: \cC \rightarrow \R$, with $\cC$ as in Definition \ref{actions}, such that 
	the functions $z_j$ depend on $a$ only, and
	$\forall a\in \Omega\cap\Lambda$ there exists
	an eigenvalue $\lambda_a$ of \eqref{ope1} which admits the
	asymptotic expansion
	\begin{equation}
	\label{asym}
	\lambda_a\sim h_L(a)+\sum_{j\geq0} z_j(a)\ .
	\end{equation}
\end{theorem}

Such a theorem ensures that in the case of multiplicity one, the
majority of the eigenvalues of the perturbed problem are stable,
namely they move only by asymptotically small quantities and
furthermore admit an asymptotic expansion.

It is also possible to give a description of all the  eigenvalues
outside $\Omega$, following the presentation of \cite{BLMres},
but the corresponding precise statement is quite complicated, so we
avoid to give it here.

\subsubsection{Time dependent perturbations}\label{growth}

The second kind of results deals with the problem of growth of Sobolev norms,
when the perturbation is time dependent. This is contained in the
following Theorem.

\begin{theorem}
	\label{main.lin}[Main Theorem of \cite{QN}]
	Let $H=H(t)$ be of the form
	\begin{equation}
	\label{H.mai}
	H(t):=H_L+V(t)
	\end{equation}
with $H_L$ fulfilling assumption \ref{L}. Assume also that $V(\cdot) \in
	C^\infty_b\left(\R; \cA^\tb \right)$ is a family of self-adjoint
pseudodifferential operators with $\tb< \td$; then for any $s\geq 0$ and for any initial datum $\psi \in \cH^s$ there exists a unique global solution $\psi(t) := \cU(t, \tau)\psi \in \cH^s$ of the initial value problem
		\begin{equation}
		\label{p.abs}
		\im \partial_t \psi(t) = H(t) \psi(t)\, ,  \quad \psi(\tau) = \psi\,,
		\end{equation}
	furthermore, for any $s>0$ and $\varep>0$ there exists a positive constant $K_{s, \varep}$ such that for any $\psi \in \cH^s$
	\begin{equation}\label{growth.eq}
	\|\cU(t, \tau) \psi\|_s \leq K_{s, \varep} \langle t -\tau \rangle^{\varep} \| \psi\|_s\ , \quad \forall t, \tau \in \R\,.
	\end{equation}
\end{theorem}

In particular Theorem \ref{main.lin} applies to the following equations:
\begin{enumerate}
\item  Quantum particle on $\T^n$ subjected to electro-magnetic potential: 
\begin{equation}\label{maxwell.edison}
\im \partial_t \psi = \sum_{j = 1}^n(\im \partial_j + B_j(t, x))^2 \psi +   W(t,x)\psi\,, \quad x \in \T^n\,,
\end{equation}
with $ B_1, \dots, B_n, W \in \cC^\infty_b(\R; \T^n)$.  This case was also studied in \cite{BLM_growth}, with techniques specific to the torus. 
\item Perturbed two-dimensional quantum anharmonic oscillator:
\begin{equation}
	\im \partial_t \psi = -\frac 1 2 \Delta \psi + \frac{1}{2\ell} \|x\|^{2\ell} \psi + \im \sum_{i, j=1,2} v_{ij}(t) x_i \partial_{x_j} \psi\,, \quad x \in \R^2
\end{equation}
with $v_{i,j}(t) \in \cC^\infty_b(\R)$.
\item Quantum particle on any one of the manifolds of
  Subsect. \ref{EGIQS} with time dependent potential:
\begin{equation}
	\im \partial_t \psi = -\Delta_{g} \psi + V(t) \psi\,, 
\end{equation}
with $V \in \cC^\infty_b(\R;M)$.
\end{enumerate}

\begin{remark}
	The classical counterpart to the problem of bounding the
        growth of Sobolev norms for \eqref{maxwell.edison} was treated
        in \cite{dario_homog} for a particle on the torus; there the following result was proven: consider the classical Hamiltonian
	\begin{equation}
		h(x, \xi) = \frac 1 2 \sum_{j=1}^n(\csi_j - B_j(t, x))^2 + W(t, x)\,, \quad (x, \xi) \in \T^n \times  \R^n\,,
	\end{equation}
	with $B_1, \dots, B_n, W \in \cC^\infty_b(\R \times \T^n)$. Then for any $\varep >0$ there exists $C_\varep >0$ such that all solutions of the Hamilton equations generated by $h$ satisfy
	$$
	\|\xi(t)\| \leq C_{\varep} \langle t \rangle^{\varep} (1 +
        \|\xi(0)\|)\ , \quad \forall t \in \R\,.
	$$
\end{remark}
\subsection{Nonlinear Hamiltonian perturbations}\label{nonlinear}

Concerning nonlinear perturbations, we have to restrict to the case of
operators on compact manifolds so we assume that $M$ is a compact
manifold, and consider a nonlinear PDE of the form
\begin{equation}
  \label{NL.eq}
\im \psi_t=H_L\psi+P(\psi,\bar\psi,x)\ ,
\end{equation}
where $P$ is a Hamiltonian vector field whose properties will be
specified in a while. Furthermore, besides assumption \ref{L}, we are
going to make a nonresonance assumption on the eigenvalues of $H_L$;
this is typically achieved by considering the case where $h_L$ depends
on some external parameters that are used to tune the
frequencies. This is the way we proceed in applications.

We start by giving the precise
assumption on $P$.

\begin{assumption}{(NL)}\label{NL}
  There exists a function $F\in C^{\infty}(\cV\times\cV\times M;\C)$ with $\cV\subset\C$ a
  neighborhood of the origin with the following properties
  \begin{itemize}
  \item[i.]  
$  P=\partial_{\bar\psi}F
    $,
  \item[ii.] there exists $C$ s.t. for any $z\in\cV$ one has $F(z,\bar z,x)\in\R$, and
    $$
|F(z,\bar z,x)|\leq C |z|^3\ .
    $$
  \end{itemize}
\end{assumption}
Actually a slightly more general assumption is enough and this is
needed for the application to the nonlinear stability of the ground
state in NLS, but we will not give here its abstract form.

We are now going to state the nonresonance assumption on the eigenvalues of
$H_L$. As we will discuss in a while, this is a generalization of
 ``0-th Melnikov condition'', whose standard form is
\begin{equation}
  \label{omel}
  \begin{gathered}
  \omega_{a^{(1)}}\pm\omega_{a^{(2)}}\pm...\pm
  \omega_{a^{(r)}}\not=0\ \quad
  \Longrightarrow  \\
  \left|\omega_{a^{(1)}}\pm\omega_{a^{(2)}}\pm...\pm
  \omega_{a^{(r)}}\right|\geq\frac{\gamma_r}{(\max|a^{(j)}|)^{\tau_r}}\ .
  \end{gathered}
\end{equation}
To give the precise statement denote by
$$
\omega_a:=h_L(a)\ ,\ a\in\Lambda
$$
the eigenvalues of $H_L$ and remark that by the assumption on the
homogeneity of $h_L$ one has a bound on the cardinality of the set of
eigenvalues in a ball. From such a bound it is easy to deduce the
following clustering property

	\begin{lemma}
		\label{parti_sigma}
		There exist two sequences $\alpha_n$, $\beta_n$ $n\in\N,
                \ n\ge0$ with
 $\alpha_n<\beta_n<\alpha_{n+1}$ and such that 
		\begin{itemize}
		\item $h_L(\Lambda)\subset \bigcup_{n}[\alpha_n,\beta_n]$;
		\item $\beta_n-\alpha_n \leq {2}$;
		\item $\alpha_{n+1}-\beta_n\geq 2 n^{-d/\td}$.
		\end{itemize}
	\end{lemma}
The nonresonance condition that we are going to assume ensures that
one can use Birkhoff normal form procedure in order to eliminate from
the nonlinearity terms
	enforcing exchanges of energy among modes labeled by indexes
	belonging to different intervals. This requires to
        characterize the linear combination of frequencies which
        correspond to monomials enforcing such exchange of
        energy. So, we 
give the following
        definition

\begin{definition}
          \label{res.index}
Given $r\in\N$ and $j\in\left\{0,...r\right\}$, a set of indexes
$(a^{(1)},...,a^{(r)})\in\Lambda^r$ will be said to be
\emph{dangerous} if $r$ is even and for any $l=1,...,r/2$ there exists
an interval $[\alpha_l,\beta_l]$ with the property that
\begin{equation}
  \label{poss}
a^{(l)} \ \text{and}\ a^{(l+r/2)}\in [\alpha_l,\beta_l]\ ,\forall
l=1,...,r/2\ .
\end{equation}
A sequence which is not dangerous is called \emph{nondangerous}. 
\end{definition}

	\begin{assumption}{(NR)}\label{NR}
		For any $r\ge3$,  there are constants $\gamma_r,\tau_r>0$ 
		such that $\forall j=1,...,r$ one has 
\begin{equation}\label{diof}
  \left|\sum_{l=1}^{j}\omega_{a^{(l)}}-\sum_{l=j+1}^{r}\omega_{a^{(l)}}\right|
  \ge \frac{\gamma_r}{\left(\max{|a^{(l)}|}\right)^{\tau_r}}\,.
		\end{equation}
for any non dangerous sequence $(a^{(1)},...,a^{(r)})$. For dangerous
sequences we require \eqref{diof} to hold only for $j\not=r/2$.
	\end{assumption}

        The following theorem is the main abstract result of \cite{QNN}.

\begin{theorem}\label{ab.res}
		Consider the Hamiltonian system \eqref{NL.eq}. Assume
		Hypotheses \ref{NL} and \ref{NR}, 
		then for any
		integer $r\geq3$, there exists $s_r\in \N$ such that, for any
		$s\ge s_r$, there are constants $\epsilon_0>0$, $c>0$ and $C$
		for which the following holds: if the initial datum
		$\psi_0\in H^s(M,\C)$ fulfills
		\begin{gather*}
			\epsilon := \left\|\psi_0\right\|_s < \epsilon_0\,,
		\end{gather*}
		then the Cauchy problem has a unique solution 
		$u\in\mathcal C^0\left((-T_\epsilon,T_\epsilon), H^s(M,\C)\right)$ with 
		$T_\epsilon > c \epsilon^{-r}$. Moreover, one has
		\begin{gather}
			\norm{\psi(t)}_s \le C \epsilon,\quad \forall
			t\in(-T_\epsilon,T_\epsilon)\ .
		\end{gather}
	\end{theorem}

	In particular, in \cite{QNN} we show that Theorem \ref{ab.res} applies to the following models:
	\begin{enumerate}
	\item Schr\"odinger equation with convolution potential on one
          of the manifolds $M$ of the examples 1, 2, 3, 4, of
          Subsection \ref{EGIQS}
	\begin{equation}\label{shro.convol}
	\im \partial_t \psi = -\Delta_g \psi + V \ast \psi + f(x, |\psi|^2) \psi\,, \quad x \in M\,,
	\end{equation}
with $f(x, \cdot)$ a smooth function having a zero of order  at least $2$  at $0$, and $V$ belonging to a full measure subset of
	$$
	\cV_m := \left\{ V = \{V_k\}_{k \in \Z^d} \ \left|\ \langle k \rangle^m |V_k| \in \left[-\frac 1 2, \frac 1 2 \right]\right.\right\} \quad \text{for some } m>0\,.
	$$
The case of $M=\T^d$ with a flat metric was already dealt with in \cite{BFM22}.
        
	\item Beam equation on one of the manifolds $M$ of the
          examples 1, 2, 3, 4, of Subsection \ref{EGIQS}
	\begin{equation}
	(\partial_{tt} + \Delta_g^2 + m )u = - \partial_u F(x, u)\,,
	\end{equation}
	with $F \in \cC^\infty(M \times \cV; \R)$, where $\cV$ is a
        neighborhood of the origin, having a zero of order at least 3
        at $u = 0$, and $m$ belonging to a full measure subset of
        $\R^+$.
        
	\item Stability of the ground state in Nonlinear
          Schr\"odinger equation on one of the above manifolds
	\begin{equation}\label{nls}
	\im \partial_t \psi = -\Delta_g \psi + f(|\psi|^2)\psi\,,
	\end{equation}
	with $f \in \cC^\infty(\cV; \R)$ having a zero of order at least 1 at the origin. Then, for any $p_0 \in \cV \cap \R^+$ the plane wave
	\begin{equation}
	\psi_\star(t) = \sqrt{p_0} e^{\im \nu t}\,, \quad \nu = f(p_0)\,,
	\end{equation}
	is a solution of \eqref{nls}. In this case, from the abstract Theorem \ref{ab.res} we deduce the following:
	\begin{theorem}
	 Let $\bar \lambda$ be the minimum nonvanishing eigenvalue of $-\Delta_g$ and assume there exists  $\bar{p}_0>0$ such that $ \bar\lambda + 2 f(p_0) > 0$
	for any $p_0 \in (0,\bar{p}_0]$.	Then there
		exists a zero measure set $\mathcal{P}\subset(0,\bar p_0]$ such that
		if $p_0\in (0,\bar{p}_0]\setminus \mathcal{P}$ then for any $r\in \N$ there
		exists $s_r$ for which the following holds. For any $s\ge s_r$, there
		exist constants $\epsilon_0$ and $C$ such that if the initial datum $\psi_0$
		fulfills
		\begin{gather*}
		\norm{\psi_0}_0^2 = p_0,\qquad \inf_{\alpha\in\T}\norm{\psi_0 - \sqrt{p_0}e^{-i\alpha}}_s = \epsilon \le \epsilon_0\,,
		\end{gather*}
		then the corresponding solution fulfills
		\begin{equation*}
		\inf_{\alpha\in\T}\norm{\psi(t) - \sqrt{p_0}e^{-i\alpha}}_s
		\le C \epsilon \qquad \forall\, |t|\le C
		\epsilon^{-r}\ .
		\end{equation*}
	\end{theorem}
	\end{enumerate}

\section{Nekhoroshev-Bourgain partition of the joint
  spectrum}\label{clustering}

In this section we give a partition of the joint spectrum $\Lambda$,
which is the quantum analogue of the partition of the action space of
a classical system in resonant blocks introduced by
Nekhoroshev in order to prove his celebrated theorem on exponential
stability of quasi integrable systems. It turns out that such a
partition also has a remarkable property that was introduced by
Bourgain in order to prove KAM type results for nonlinear wave and
Schr\"odinger equation on $\T^d$.

As a result, this partition is the main tool for proving all the results
of the previous section.

\subsection{Classical Nekhoroshev Theorem and classical Nekhoroshev's
  Partition} \label{le.ode.sono.utili} 
We start with recalling the classical Nekhoroshev Theorem:
\begin{theorem}\label{nek.classico}
	Let $h: \cA \times \T^n \rightarrow \R$ be a function of the form
	$$
	h(a, \alpha) = h_L(a) + \varep v(a, \alpha)\,,
	$$ where $h_L$ and $v$ are analytic on $\cA$ and $h_L$ is
        steep on $\cA$.
        Then there exist
        $\varep_0>0$ $\mathtt{a}, \tb>0$  such that, defining $\cU := \{a \in \cA \ |\ \operatorname{dist}(a, \partial \cA) \geq 2 \varep^\mathtt{b}\}$, if $|\varep|
        < \varep_0$, then for all initial data in $\cU\times\T^n$ one
        has
	\begin{equation}
	\|a(t) - a(0)\| \leq \varep^{\mathtt{b}} \quad \forall |t| \leq \varep^{-1 }e^{ \varep^{-{\mathtt{a}}}}\,.
	\end{equation}
	
\end{theorem}

The classical proof as formulated by Nekhoroshev in \cite{Nek77, Nek79} (see also \cite{GCB, gio_pisa, nek_noi, dario_homog}) is based on two steps: an analytic part, in which the Hamiltonian $h$ is put in local normal form, removing all non resonant contributions, and a geometric part, based on the analysis of resonances in the space of actions.
 We now are going to describe the latter geometric analysis of resonances and its role in the proof of Nekhoroshev Theorem, in order to enlighten the similarities with the procedure we adopt in the quantum case.

Actually, here we present the scheme of the proof of
\cite{nek_noi,tesi_silvia}, which is the closest to the quantum case
and deals with long time stability in the case of smooth Hamiltonians:
\begin{theorem}[$\cC^\infty$ Nekhoroshev Theorem]\label{nek.liscio}
	Let $h:\cA \times \T^n \rightarrow \R$ be of the form
	$$
	h(a, \alpha) = h_L(a) + \varep v(a, \alpha)\,, 
	$$
	with $h_L$ steep on $\cA$ and $h_L,v \in \cC^\infty_b(\cA \times \T^n; \R)$. For any
        $b < \frac 1 2$ and $N>0$ there exist $\varep_N>0$ and $C_N>0$
        such that, for $0< \varep < \varep_N$, and any initial datum
        in $\cU\times \T^n$, with $\cU := \{a \in \cA \ |\ \operatorname{dist}(a, \partial \cA) \geq 2 \varep^{b}\}$, one has
	\begin{equation}
		\|a(t) - a(0)\| \leq \varep^{b} \quad \forall |t| \leq C_N \varep^{-N}\,.
	\end{equation}
\end{theorem}

The first step in the proof of Theorem \ref{nek.liscio} consists in performing a normal form procedure, which for any $N$ enables to conjugated the Hamiltonian $h$  to a new one of the form
$$
h_N = h_L(a) + \varep z_N(a,\alpha) + \varep^N r_N(a, \alpha)\,,
$$
with
\begin{equation}\label{su.tutte}
	z_N(a, \alpha) = \sum_{k \in \Z^n \atop \|k\| \leq K} \hat{z}_k(a) e^{\im k \cdot \alpha}\,,
\end{equation}
$K= \varep^{-\frac{1}{\zeta}}$ for some $\zeta >0$. Furthermore, there exists $\delta >0$ such that, for each $k$,
\begin{equation}\label{le.fourier}
	\operatorname{supp} \hat{z}_k(a) \subseteq \{ a \ |\ |w(a) \cdot k| \leq \varep^{\delta} \}\,,
\end{equation}
where we recall that $w(\cdot) := \partial_a h_L(\cdot)$. Then the
geometric part of the proof consists in analyzing the dynamics
generated by $h_L + \varep z_N$.

To this aim we first recall that a module of dimension $s=0, \dots,
n$, is a subset $M \subseteq \Z^n$ such that
$\operatorname{span}_\R(M) \cap \Z^n = M$. In this section we will
denote by $M$ a submodulus of $\Z^n$.

Now we start by defining the \emph{resonant zones}
\begin{multline}
	\cZ^{(s)}_{M} := \{ a \ |\ \exists\  k_1, \dots, k_s \quad \text{linearly indep. in }M \\
	\text{s.t.} \quad  |w(a) \cdot k_j| \leq \varep^{\delta_s} \quad \forall j = 1, \dots, s \}\,,                      
\end{multline}
with a suitable sequences of increasing $\{\delta_s\}_{s=1}^{n}$ with $\delta_1 = \delta$.
For any $M$, the set $\cZ^{(s)}_M$ collects all points $a$ such that $w(a)$ is resonant with some vectors in $M$. Its role is that, due to the form of the truncated Hamiltonian $h_L + \varep z_N$ (see \eqref{su.tutte}--\eqref{le.fourier}), points in $\cZ^{(s)}_M$ only move along directions parallel to $M$. 
Then, as long as the point $a(t)$ belongs to some resonant zone, its motion is easy to describe. A key step in the geometric part of the proof is then to understand how the resonant zones intersect one each other, and whether or not is possible to pass from one resonant zone from the other, visiting sooner or later all the action space (this is the well known phenomenon of \textit{overlapping of resonances}).

	\begin{figure}[H]
	\begin{subfigure}[t]{0.47\textwidth}
		\includegraphics[width=\textwidth]{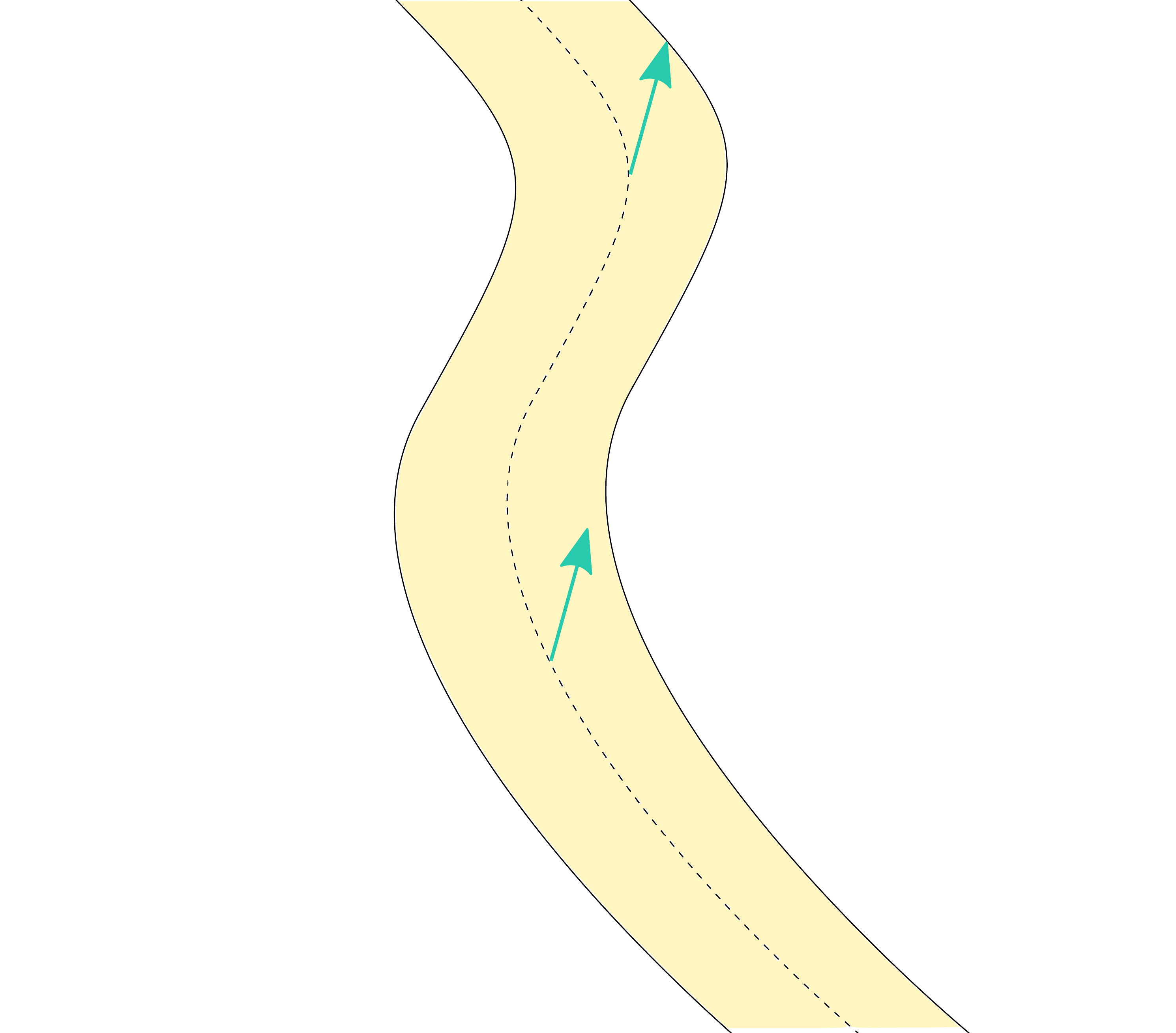}
		\caption{In yellow, cartoon of a resonant zone $\cZ^{(s)}_M$, with $M=\operatorname{span}\{k\}$ for a given $k$. The  dashed line corresponds to the curve $w (a)\cdot k =0$. Inside the zone, the motion generated by $h_L + \varep z_N$ is parallel to $k$, pictured in turquoise. }
		\label{1.zone}
	\end{subfigure}
	\hfill
	\begin{subfigure}[t]{0.47 \textwidth}
		\includegraphics[width=\textwidth]{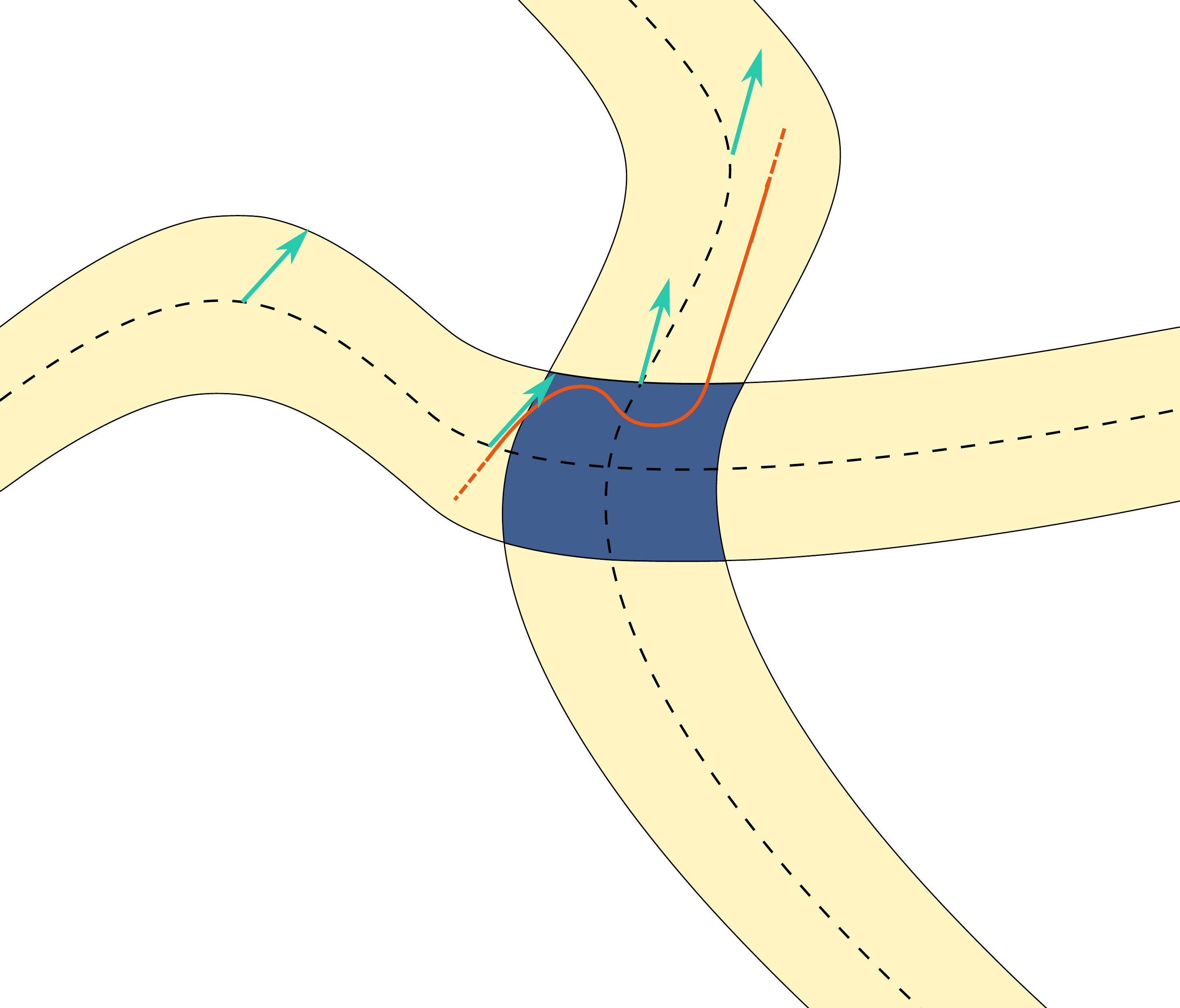}
		\caption{In yellow, two resonant zones associated to different modules. Their intersection is pictured in blue. In orange, a possible dangerous trajectory, passing from one zone to the other (overlapping of resonances, which we want to exclude).
		}
		\label{blocco}
	\end{subfigure}
\end{figure}

In order to show that the latter phenomenon does not occur, one needs to define a more refined decomposition of the actions space $\cA$ with respect to the one given by the resonant zones $\{\cZ^{(s)}_{M}\}_{s, M}$. In particular, one proves the existence of a partition
\begin{equation}
	\cA = \bigcup_{s=0}^{n} \bigcup_{\begin{subarray}{c}		M \subseteq \Z^n \\ \operatorname{dim} M = s \\ M \text{module}
		\end{subarray}} \bigcup_{j \in \N} \mathcal{E}^{(s)}_{M, j}
\end{equation}
		in sets $\cE^{(s)}_{M, j}$ such that
		\begin{enumerate}
		\item Each $\cE^{(s)}_{M, j}$ has small diameter:
		\begin{equation}
		\operatorname{diam}\left(\mathcal{E}^{(s)}_{M, j}\right) \leq \varep^{\mathtt{c}_2} \quad \forall s, M, j\,.
		\end{equation}
		\item Each $\cE^{(s)}_{M, j}$ is contained in $\cZ^{(s)}_{M}$, and is far away from having resonances with vectors not in $M$: in particular, if $k \in \Z^n$ satisfies
		\begin{equation}
		|w (a) \cdot k| < \varep^{\delta_1}\,, \quad \|k\| \leq K\,,
		\end{equation}
		then $k \in M$ and $a + \eta k \in \cE^{(s)}_{M, j}$,
                provided $|\eta | \ll 1$.
               
           \item $\cE^{(s)}_{M, j}$ is invariant for the dynamics
                  of $h_L+\varepsilon z_N$.
 
		\end{enumerate}
	Note that, once this decomposition has been obtained, the statement of Theorem \ref{nek.liscio} essentially follows combining Items 1 and 3.
	
\subsection{Invariant partition in the quantum case}
Let us start with defining what we consider \emph{resonant} in the case of a quasi-integrable quantum system:
\begin{definition}\label{def.res}
	Let $\delta < \td - 1$, $\mu >0$ and $\tR >0$. We say $a\in \Lambda$ \emph{is resonant with} $k \in \Z^d$ if
	\begin{equation}\label{risuono.QN}
	|w (a) \cdot k| < \|a\|^{\delta} \|k\|\,, \quad \|k\| \leq \|a\|^\mu\,, \quad \|a\| \geq \mathtt{R}\,,
	\end{equation}
        where $w(a):=\partial_ah_L(a)$
\end{definition}
A key point in proving Theorems \ref{maint}, \ref{rilinea} and
\ref{main.lin} is the following result, proven in \cite{QN} which
gives the partition constituting the main step for the proof of the
theorems of Sect. \ref{perturbations}:

\begin{theorem}\label{teo.partizione}
	Under the assumptions of Theorem \ref{main.lin}, there exist $0<\delta_0< \td-1$, $\mu_0>0$ and $\mathtt{R}_0>0$ such that, if $\delta_0<\delta < \td - 1$, $0<\mu<\mu_0$, and $\mathtt{R}>\mathtt{R}_0$, the lattice $\Lambda$ admits a partition
	\begin{equation}
	\Lambda = \bigcup_{s=0}^{d} \bigcup_{\begin{subarray}{c} 	M \subseteq \Z^d \\ \operatorname{dim} M = s \\ M \text{module}
		\end{subarray}} \bigcup_{j \in \N} E^{(s)}_{M, j}
	\end{equation}
	with the following properties:
	\begin{itemize}
	\item[(i)] Each block $E^{(s)}_{M,j}$ is dyadic, in the sense that there exists $C>0$ such that
	\begin{equation}\label{dyadic}
	 \max\{ \| a \| \ |\ a \in E^{(s)}_{M,j}\} \leq C \min \{ \| a \| \ |\ a \in E^{(s)}_{M,j}\}
	\end{equation}
	\item[(ii)]  If $a \in E^{(s)}_{M,j}$, there exists a vector $k \in \Z^d$ such that $a$ is resonant with $k$. Conversely, if $a$ \emph{is resonant with} $k \in \Z^d$ in the sense of Definition \ref{def.res}, 
	then $a + k \in E^{(s)}_{M,j}$.
	\item[(iii)] There exists $K>0$ such that, if $a\in E^{(s)}_{M, j}$ and $b\in E^{(s')}_{M', j'}$ with $(s, M, j) \neq (s', M', j')$, then
	\begin{gather}\label{lontano.lontano}
	\|a-b\| + |\omega_a-\omega_b| \ge K (\|a\|^\mu + \|b\|^\mu)\,.
	\end{gather}
	\end{itemize}
\end{theorem}

We start by remarking the analogies with the properties of the
covering $\{\cE^{(s)}_{M, j}\}$ exhibited in the classical case in
Subsection \ref{le.ode.sono.utili}: Item $(i)$ gives an upper bound on
the diameter of the $E^{(s)}_{M,j}$, just as Property 1 of the sets
$\cE^{(s)}_{M,j}$. Furthermore, Item $(ii)$ states that, if $a$ is
resonant with $k$, and $\|k\|$ is not too large, then $a + k$ belongs
to the same block $E^{(s)}_{M,j}$, which is analogous to Property 2 of
the sets $\cE^{(s)}_{M,j}$.

We also remark that a first connection between Nekhoroshev's partition
and the spectral properties of the perturbations of the Laplacian was
pointed out in \cite{FFS}.

\vskip10pt

{\it         We now describe how Theorem \ref{teo.partizione} is used to
        prove the different theorems.
}
        \subsection{Linear perturbations}
        The first step of the proof consists in  proving a
        ``quantum normal form theorem'' for Hamiltonians of the for
        \eqref{ope1} and \eqref{H.mai}.

        		This is contained in the following:
		\begin{proposition}\label{forma.formale}
		For any $N \in \N$ there exists a map $U_N \in \cC^\infty_b(\R; \cB(\cH^s, \cH^s))$, such that $U_N(t)$ is unitary in $\cH$ for any $t$ and $\psi$ solves \eqref{p.abs} if and only if $\phi = U^{-1}_N(t) \psi$ solves
		\begin{equation}\label{fn.1}
		\im \partial_t \phi = (H_L + Z_N(t) + R_N(t)) \phi\,,
		\end{equation}
		with $Z_N$ and $R_N$ self-adjoint operators as follows:
		\begin{enumerate}
		\item \emph{(Smoothing remainder)} $R_N \in \cC^\infty_b(\R; \cB(\cH^s; \cH^{s + N}))$ for any $s$
		\item \emph{(Resonant normal form)}  $Z_N \in \cC^\infty_b(\R; \cB(\cH^s; \cH^{s -\tb}))$ for any $s$. Moreover,  for any $a \in \Lambda$ let $\Pi_a$ be the projection operator defined in \eqref{Pi.a}, then one has
		\begin{equation}\label{matr.el}
		\Pi_{a} Z_N(t) \Pi_{a+k} \neq 0
		\end{equation}
		only if $a$ is resonant with $k$, or $a + k$ is resonant with $k$.
		\item \emph{(Block diagonal structure)} As a consequence, $Z_N$ is block diagonal along the blocks $E^{(s)}_{M,j}$ of Theorem \ref{teo.partizione}, in the sense that for any $t$
		\begin{equation}
			\Pi_{a} Z_N(t) \Pi_b \neq 0 \quad \Rightarrow a, b \in E^{(s)}_{M,j}
		\end{equation}
		for some $s, M, j$.
		\end{enumerate} 
		\end{proposition}
		Proposition \ref{forma.formale} is based on an iterative procedure along which, for any $n$, at the $n-$th step  one is able to eliminate from the perturbation $V_n(t)$ all terms
		$$
		\Pi_{a} V_n(t) \Pi_{a + k}
		$$
      	such that both $a$ and $a + k$ are \emph{not resonant} with $k$. Then all the remaining resonant contributions are collected in the normal form $Z_N$. This is the heart to the proof of Item 2. {Note also that Item 3. follows directly combining Item 2. with property $(ii)$ of Theorem \ref{teo.partizione}. }
      

	We now discuss how Proposition \ref{forma.formale} enables to Prove Theorems \ref{maint}, \ref{rilinea} and \ref{main.lin}.
\begin{proof}[Sketch of the proof of Theorems \ref{rilinea} and
    \ref{maint}] We apply Proposition \ref{forma.formale} in the time independent case. The main point is to define $\Omega$ as the set
  $$
\Omega:=\bigcup_{j}E^{(0)}_{j,M}
  $$
where $M=\left\{0\right\}$ is the trivial modulus. By Item $(ii)$ of Theorem \ref{teo.partizione}, it follows that such a set
contains all the points $a$ with the property that $w(a)$ is non-
resonant, namely the only vector $k\in \Z^d$ such that \eqref{risuono.QN} is satisfied is $k=0$. As a consequence, by Item 2 of Proposition \ref{forma.formale} at these points the normal form $Z_N$ turns
out to leave invariant the common eigenspaces of the actions,
namely one has
$$
Z_N\Pi_a=\Pi_aZ_N\Pi_a\ ,\quad \forall a\in\Omega\ .
$$
Furthermore the algorithm of the proof of Proposition \ref{forma.formale} also
allows to show that $Z_N=\langle V\rangle+$lower order operators, with $\langle V \rangle$ defined by \eqref{average}.
This is why the spectrum of the restriction of $\langle V\rangle$ to
$\Pi_a\cH$ controls the first correction to the unperturbed
eigenvalues. Then one should prove that the set $\Omega$ has density
one at infinity. This can be obtained by following some ideas of the
proof of the corresponding theorems of \cite{BLMnr} and \cite{anarmonico}, but exploiting the steepness of $h_L$ instead than the arguments in \cite{Rus01}.

The proof of Theorem \ref{maint} then simply follows by the remark
that in this case the pseudo-differential operators commuting with all
the actions are just functions of the actions.
\end{proof}

\begin{proof}[Sketch of the proof of Theorem \ref{main.lin}]

  By Item 3. of Proposition \ref{forma.formale}, neglecting the
  remainder $R_N(t)$, one can study the dynamics generated by $H_L +
  Z_N(t)$: by the fact that (1) $Z_N$ is diagonal with respect to the
  blocks, (2) $Z_N$ is self-adjoint so that the $L^2$ norm is conserved
  along the dynamics it generates, (3) since the blocks are dyadic, in
  each block the $L^2$ norm is equivalent to the $\cH^s$ norm, one
  gets that, along the dynamics of $H_L +
  Z_N(t)$ one has
  	$$
	\|  \varphi(t)\|_{\cH^s} \lesssim_s  \| \varphi(0)\|_{\cH^s}\,.
	$$
	Then the contribution of the smoothing remainder $R_N(t)$ is
        dealt with using Duhamel principle. Precisely, denote by $\cU$
        the flow map of the complete Hamiltonian \eqref{fn.1} and by
        $\cU_{app}$ the flow map of $H_L +
  Z_N(t)$, then,
        by Duhamel formula
        one has
	\begin{equation}
	\cU(t, \tau) \psi = \cU_{app}(t, \tau) \psi + \int_{\tau}^{t}
	\cU_{app}(t, \tau^\prime) (-\im R^{(N)}(\tau^\prime))
	\cU(\tau^\prime, \tau) \psi\ d \tau^\prime\,,
	\end{equation}
	Since $R^{(N)}\in\timereg\left(\R; \cB(\cH^{0};\cH^{ N})\right)$ and since $H^{(N)}$ is self-adjoint in $L^2$, 
	one gets 
	\begin{align*}
	&\|\cU(t,\tau) \psi\|_{ N} \sleq
 \| \psi \|_{ N} +
	\int_{\tau}^{t} 
	\left\|R^{(N)}\right\|_{\cB(\cH^0;\cH^{ N})} \| \cU_{app}(\tau^\prime, \tau)\psi\|_{0}\ d \tau^\prime\\
	&\sleq  \| \psi \|_{ N} +
	\int_{\tau}^{t}
	\left\|R^{(N)}\right\|_{\cB(\cH^0;\cH^{
			N})} 
	\| \psi\|_{0}\ d \tau^\prime\sleq  \| \psi \|_{ N} +
 \langle t -\tau\rangle \| \psi \|_0 \,,
 \\
& \sleq \langle t -\tau\rangle \| \psi \|_N\ .
        \end{align*}
Then using interpolation one
        immediately gets the thesis.
      \end{proof}
  
\subsection{Nonlinear perturbations}

\begin{proof}[Sketch of the proof of Theorem \ref{ab.res}]
Theorem \ref{ab.res} is the consequence of a Birkhoff normal form
theorem conjugating \eqref{NL.eq} to its linear part plus a ``normal
form'' plus a remainder which is of order $\epsilon^N$, with
arbitrary $N$. Now, the main point is that, as it is well known,  in
order to perform such a normal form one has to consider small
divisor corresponding to the terms one wants to eliminate from
the Hamiltonian. The kind of
normal form one can obtain is thus directly related to the kind of
nonresonant conditions which is satisfied by the frequencies
$\omega_a$. Using just the nonresonant condition \eqref{omel} it is
difficult to obtain a useful normal form in the general
case. Actually the kind of nonresonant condition which is known to be
very useful is 
\begin{equation}
  \label{omel.1}
  \begin{gathered}
 \omega_{a^{(1)}}\pm\omega_{a^{(2)}}\pm...\pm
 \omega_{a^{(r)}}\not=0\ \Longrightarrow \\
 \left|\omega_{a^{(1)}}\pm\omega_{a^{(2)}}\pm...\pm
 \omega_{a^{(r)}}\right|\geq\frac{\gamma_r}{(\max_3|a^{(j)}|)^{\tau_r}}\,,
  \end{gathered}
\end{equation}
where $\max_3$ is the third largest index among the $a^{(j)}$'s.

The idea here is that if one only tries to eliminate from the
nonlinearities terms coupling different blocks $E_{M,j}$, then it
turns out that for the corresponding small denominators one has a
selection rule, which entails the fact that condition \eqref{omel}
implies condition \eqref{omel.1} for all the small denominators that
really matter.

The actual implementation of such a scheme is however quite
complicated technically and we refer to \cite{QNN} for more details
(see also \cite{BFM22}).
\end{proof}

\section{Two examples of quantum integrable systems}\label{sez.per.esempio}
In this section we analyze in detail two models and we show that they satisfy the assumptions of globally integrable quantum systems.
\subsection{Anharmonic oscillator}
Let us consider the operator
\begin{equation}
H_{An} = -\frac{\Delta}{2} + \frac{\left\|x\right\|^{2\ell}}{2\ell}\,, \quad x \in \R^2\,,
\end{equation}
with $\ell \in \N$, $\ell >1$. 
In this case,  the class $\cA^m$ of pseudodifferential operators that we consider is the one defined in the second line of  \eqref{psido}, and we are going to show that there exist a function $h_L$ and two action operators $A_1, A_2 \in \cA^1$ such that
$$
H_{An} = h_L(A_1, A_2) + R\,,
$$
where $h_L$ coincides outside a neighborhood of the origin with the
classical Hamiltonian written in action variables and $R$ is a lower
order correction. It was proved in
\cite{QN} that it is homogeneous of degree $\td = \frac{2\ell}{\ell +
  1}$; the operator $R$ is a lower order operator in $\cA^\tb$ with
$\tb < \td$. Remark that lower order terms are not relevant in order
to satisfy the assumptions of Theorem \ref{main.lin}, since the term
$R$ can be included in the perturbation $V(t)$.

To prove this fact, following the strategy of Section \ref{flows}: we
start by remarking that $H_L$ is the quantization of the classical Hamiltonian
\begin{equation}
\label{ana.classica}
 h_{An}(x,\xi)=\frac{\left\|\xi\right\|^2}{2}+
\frac{\left\|x\right\|^{2\ell}}{2\ell} \,,
\end{equation}
and start by constructing suitable action angle variables for it. The action variables that we choose for $h_{An}$ are the angular momentum
\begin{equation}
\label{az.L}
a_2(x,\xi):=x_1\xi_2-x_2\xi_1\ 
\end{equation}
and the radial action $a_r$, namely the action of the effective Hamiltonian
\begin{equation}
\label{efficace}
h_{An}^*(r,p_r, \tM):=\frac{p_r^2}{2}+V_{\tM}^*(r)\ ,\quad
V_{\tM}^*(r):=\frac{\tM^2}{2r^2}+\frac{r^{2\ell}}{2\ell}\ ,
\end{equation}
where $r={\|x\|}$, $p_r$ is the conjugated
momentum and $\tM$ is a particular value of the angular momentum, in the sense that we put ourselves on a level surface {$\{(x, \xi)\ |\ a_2(x,\xi)=\tM\}$}.
In particular, one has
\begin{equation}
\label{ar}
a_r=a_r(E,\tM):=\frac{\sqrt2}{\pi} \int_{r_m}^{r_M}\sqrt{ E- V^*_\tM(r)}dr\ ,
\end{equation}
where $0<r_m <r_M$ are the two solutions of the equation
$$
E-V^*_\tM(r)=0\ .
$$ However, the actions $(a_2, a_r)$ are \emph{not} good candidates to
be quantized for our construction: indeed, while $a_2$ is a
polynomial, $a_r$ turns out to have a singularity at $\tM=0$.  Anyway,
this singularity can be removed, by performing a suitable unimodular
transformation. This is the content of the following lemma, proved in
\cite{anarmonico}:
\begin{lemma}
	\label{azioni.osci}[Lemma 4.5 of \cite{anarmonico}] 
	The function
	\begin{equation}
	\label{a1}
	a_1(E,\tM):=\left\{
	\begin{matrix}
	a_r(E,\tM) & for & \tM>0
	\\
	a_r(E,\tM)-\tM & for & \tM< 0
	\end{matrix}
	\right.   
	\end{equation}
	has the following properties:
	\begin{itemize}
		\item[(1)]  it extends to a complex analytic function of $\tM$ and $E$
		in the region
		\begin{equation}
		\label{defia1.1}
		\left|\tM\right|<\left(\frac{2\ell}{\ell+1}E\right)^{\frac{\ell+1}{2\ell}}\ ,\quad
		E>0\ .
		\end{equation}
		\item[(2)] Define the cone $\cC$ by
		\begin{equation}
		\label{pigreco}
		\cC:=\left\{a\in\R^2\ ;\ a_1\geq 0\ \textrm{ if } a_2 \geq0\,,\quad  a_1 \geq |a_2|\ \textrm{ if } \ a_2< 0 \right\}\,.
		\end{equation}
		Then $h_{An}$ is an analytic function of $(a_1, a_2)$ in the interior of $\cC$. Furthermore it is homogeneous of degree
		$\frac{2\ell}{\ell+1}$ as a function of $(a_1,a_2)$.
		\item[(3)] Let $a = (a_1, a_2)$, there exist  positive constants $C_1,C_2$ s.t.
		$$
		C_1\langle a\rangle \leq \tk_0\leq C_2 \langle a\rangle \ .
		$$
	\end{itemize}
\end{lemma}
Still, both $h_{An}$ and $a_1$ have a singularity at the origin. But, essentially due to the fact that the behavior of functions in a finite neighborhood of the origin does not play any fundamental role in symbolic calculus, this is easily fixed: we take $K\gg 1$ and we consider the decomposition
		\begin{equation}
		\label{regola}
		h_{An}(x,\xi)=h_{An}(x,\xi)\chi(K h_{An}(x,\xi))+h_{An}(x,\xi)(1-\chi(Kh_{An}(x,\xi)))\,.
		\end{equation}
		Now the first term is smoothing and, denoting
		\begin{equation}\label{h0reg}
		h_L(x,\xi):=h_{An}(x,\xi)(1-\chi(Kh_{An}(x,\xi)))\,,
		\end{equation}
actually $h_L$ is the expression of the function \eqref{h0reg} in action variables
$a'_1, a_2$, so that it coincides with $h_{An}$ (when written in terms of
action variables) outside a neighborhood of the origin. Moreover, since by Item $(2)$ of Theorem \ref{azioni.osci} $h_{An}$ is a homogeneous function of $(a_1, a_2)$, $h_L$ is homogeneous at infinity.
		
		Turning back to the quantum case, Theorem 3.1 and Theorem 3.2 of \cite{CdV}  guarantee the following:

\begin{theorem}
	\label{3.1cdv} There exist two pseudodifferential operators $R_j \in \cA^{\rho}$, with $\rho<1$,  and two commuting pseudodifferential
	operators $A_j\in \cA^1$ such that $A_j = Op^W(a_j) + R_j $, and $A_1, A_2$ are quantum actions according to Definition \ref{actions}. In
	particular 
	there exists a symbol $h\in S^{\frac{2\ell}{\ell+1}}_{AN}$ s.t.
	\begin{equation}\label{ham.laplaciano.1}
	H_{An}=h(A)\ .
	\end{equation}
	Furthermore one has an asymptotic expansion
	$$
	h=h_L+ \text{lower order terms}.
	$$
	with $h_L$ as in \ref{h0reg}.
\end{theorem}
 It remains to prove that $h_L$ satisfies Assumption L. This is done
 in the following steps: first one proves that there exists an open
 cone $\cC_e$ containing $\overline{\cC} \setminus \{0\}$ such that
 $h_L$ admits an homogeneous at infinity and analytic extension on
 $\cC_e$. Then one checks that $h_L$ is steep on the cone $\cC_e$;
 this is done in \cite{QN} by applying the criterion given in Theorem
 \ref{nie06}. This is verified following the methods of \cite{BFS},
 which consist in analyzing the expansion of the Hamiltonian at the
 circular orbits.
 In particular, we prove that the restriction of $h_L$ on any affine line only admits isolated critical points.

\subsection{Lie Groups}
	Consider the Laplace Beltrami operator $-\Delta_g$ in $(M, g)$, where $M$ is a simply connected Lie group endowed with the bi-invariant metric $g$. Here the class of pseudo-differential operators $\cA^{m}$ is given by the first line of \ref{psido}, namely H\"ormander class. The discussion in Section \ref{Lie} shows that, if $A_1, \dots, A_{d}$ are operators defined in \eqref{az.Lie} and $h_L : \R^d \rightarrow \R$ is the function defined by
	$$
	h_L(a) := \sum_{i, j = 1}^{d} a_i a_j \tf_{i} \cdot \tf_{j}\,,
	$$
	with $\tf_1, \dots, \tf_d$ the basis of fundamental weights of $M$ (see \eqref{cono.pesi}), then one has
	\begin{equation}
	-\Delta_g = h_L(a) - F\,, \quad F := \sum_{i, j = 1}^{d} \tf_i \cdot \tf_j\,.
	\end{equation}
	Now, as pointed out in Remark \ref{Lie.2}, the results of \cite{Fischer} prove that $A_1, \dots, A_d$ are pseudo-differential operators in $\cA^1$. Moreover, since they act as Fourier multipliers on the basis of the fundamental weights $\{f_j\}_{j}$, they mutually commute and their joint spectrum is given by elements $\Lambda = \Lambda_+(M)$ in the cone of dominant weights of $M$.
	Since one also has
	$$
	-\Delta_g  \leq -\Delta_g + |F| \leq |h_L(A)| + 2 |F| \leq \sum_{j=1}^d A^2_j \|\tf_j\|^2 + 2 |F| \lesssim \uno + \sum_{j=1}^d A^2_j\,,
	$$
	we deduce that $\uno + \sum_{j=1}^d A^2_j \lesssim -\Delta_g = K_0^2$. Hence we can conclude that $A_1, \dots, A_d$ are a set of quantum actions according to Definition \ref{actions}. We then have that $H_L = h_L(A) = -\Delta_g + F$ is globally integrable. Indeed, the function $h_L$ is homogeneous of degree two and convex, thus in particular it is steep and satisfies Assumption \ref{L}.


\begin{thebibliography}{10}
	
	\bibitem[Bam03]{Bam03}
	D.~Bambusi, \emph{Birkhoff normal form for some nonlinear {PDE}s}, Comm. Math.
	Phys. \textbf{234} (2003), no.~2, 253--285. \MR{1962462}
	
	\bibitem[Bam17]{Bam2}
	\bysame, \emph{Reducibility of 1-d {S}chr\"odinger equation with time
		quasiperiodic unbounded perturbations, {II}}, Comm. Math. Phys. \textbf{353}
	(2017), no.~1, 353--378. \MR{3638317}
	
	\bibitem[Bam18]{Bam1}
	\bysame, \emph{Reducibility of 1-d {S}chr\"odinger equation with time
		quasiperiodic unbounded perturbations, {I}}, Trans. Amer. Math. Soc.
	\textbf{370} (2018), no.~3, 1823--1865.
	
	\bibitem[Bam24]{dario_homog}
	\bysame, \emph{A global {N}ekhoroshev theorem for particles on the torus with
		time dependent {H}amiltonian}, arXiv preprint: 2401.02822, 2024.
	
	\bibitem[BD18]{BD}
	M.~Berti and J.M. Delort, \emph{{A}lmost {G}lobal {S}olutions of
		{C}apillary-gravity {W}ater {W}aves {E}quations on the {C}ircle}, UMI Lecture
	Notes, Springer, 2018.
	
	\bibitem[BDGS05]{BDGS}
	D.~Bambusi, J.M. Delort, B.~Gr{\'e}bert, and J.~Szeftel, \emph{Almost global
		existence for {H}amiltonian semilinear {K}lein‐{G}ordon equations with
		small {C}auchy data on {Z}oll manifolds}, Communications on Pure and Applied
	Mathematics \textbf{60} (2005).
	
	\bibitem[Bea77]{beals}
	R.~Beals, \emph{Characterization of pseudodifferential operators and
		applications}, Duke Math. J. \textbf{44} (1977), no.~1, 45--57. \MR{435933}
	
	\bibitem[Bel85]{Bel}
	J.~Bellissard, \emph{Stability and instability in quantum mechanics}, Trends
	and developments in the eighties ({B}ielefeld, 1982/1983), World Sci.
	Publishing, Singapore, 1985, pp.~1--106. \MR{853743}
	
	\bibitem[BFLM24]{QNN}
	D.~Bambusi, R.~Feola, B.~Langella, and F.~Monzani, \emph{Almost global
		existence for some {H}amiltonian {PDE}s on manifolds with globally integrable
		geodesic flow}, arXiv preprint: 2402.00521, 2024.
	
	\bibitem[BFM24]{BFM22}
	D.~Bambusi, R.~Feola, and R.~Montalto, \emph{Almost global existence for some
		{H}amiltonian {PDE}s with small {C}auchy data on general tori},
	Communications in Mathematical Physics \textbf{405} (2024), no.~15, 253--285.
	
	\bibitem[BFS18]{BFS}
	D. Bambusi, A. Fus\`e, and M. Sansottera, \emph{Exponential
		stability in the perturbed central force problem}, Regul. Chaotic Dyn.
	\textbf{23} (2018), no.~7-8, 821--841. \MR{3910168}
	
	\bibitem[BG01]{BG01}
	D. Bambusi and S. Graffi, \emph{Time quasi-periodic unbounded
		perturbations of {S}chr\"{o}dinger operators and {KAM} methods}, Comm. Math.
	Phys. \textbf{219} (2001), no.~2, 465--480. \MR{1833810}
	
	\bibitem[BG03]{BG03}
	D. Bambusi and B. Gr\'{e}bert, \emph{Forme normale pour {NLS} en
		dimension quelconque}, C. R. Math. Acad. Sci. Paris \textbf{337} (2003),
	no.~6, 409--414. \MR{2015085}
	
	\bibitem[BG06]{BG06}
	D.~Bambusi and B.~Gr{\'e}bert, \emph{Birkhoff normal form for partial
		differential equations with tame modulus}, Duke Mathematical Journal
	\textbf{135} (2006), no.~3, 507 -- 567.
	
	\bibitem[BGMR18]{BGMR1}
	D.~{Bambusi}, B.~{Grebert}, A.~{Maspero}, and D.~{Robert}, \emph{{Reducibility
			of the quantum Harmonic oscillator in $d$-dimensions with polynomial time
			dependent perturbation}}, Analysis \& {PDE}s \textbf{11} (2018), no.~3,
	775--799.
	
	\bibitem[BGMR21]{BGMR2}
	D. Bambusi, B. Gr\'{e}bert, A. Maspero, and D. Robert,
	\emph{Growth of {S}obolev norms for abstract linear {S}chr\"{o}dinger
		equations}, J. Eur. Math. Soc. (JEMS) \textbf{23} (2021), no.~2, 557--583.
	\MR{4195741}
	
	\bibitem[BL20]{nek_noi}
	D. Bambusi and B. Langella, \emph{A ${C}^\infty$ {N}ekhoroshev
		theorem}, Mathematics in Engineering \textbf{3} (2020), no.~2, 117.
	
	\bibitem[BL22]{QN}
	D.~Bambusi and B.~Langella, \emph{Growth of {S}obolev norms in quasi integrable
		quantum systems}, arXiv preprint: 2202.04505, 2022.
	
	\bibitem[BLM20]{BLMnr}
	D.~Bambusi, B.~Langella, and R.~Montalto, \emph{On the spectrum of the
		{S}chr{\"o}dinger operator on ${T}^d$: a normal form approach},
	Communications in Partial Differential Equations \textbf{45} (2020), 1--18.
	
	\bibitem[BLM22a]{BLM_growth}
	\bysame, \emph{Growth of {S}obolev norms for unbounded perturbations of the
		{S}chr\"{o}dinger equation on flat tori}, J. Differential Equations
	\textbf{318} (2022), 344--358. \MR{4387286}
	
	\bibitem[BLM22b]{BLMres}
	\bysame, \emph{{Spectral asymptotics of all the eigenvalues of Schr{\"o}dinger
			operators on flat tori}}, Nonlinear Analysis \textbf{216} (2022), 112679.
	
	\bibitem[BLR22]{anarmonico}
	D.~Bambusi, B.~Langella, and M.~Rouveyrol, \emph{On the stable eigenvalues of
		perturbed anharmonic oscillators in dimension two}, Communications in
	Mathematical Physics \textbf{390} (2022), no.~1, 309--348.
	
	\bibitem[BM18]{Bam3}
	D.~Bambusi and R.~Montalto, \emph{Reducibility of 1-d {S}chr\"odinger equation
		with time quasiperiodic unbounded perturbations, {III}}, J. Math. Phys.
	\textbf{59} (2018), no.~12, 122702.
	
	\bibitem[BMM22]{BMM}
	M.~Berti, A.~Maspero, and F.~Murgante, \emph{Hamiltonian birkhoff normal form
		for gravity-capillary water waves with constant vorticity: almost global
		existence}, arXiv preprint:2212.12255, 2022.
	
	\bibitem[Cai21]{tesi_silvia}
	S.~Caimi, \emph{The ${C}^\infty$ {N}ekhoroshev theorem}, Master Thesis,
	{U}niversita' degli {S}tudi di {M}ilano, 2021.
	
	\bibitem[CDV]{CdVBook}
	Y.~Colin De~Verdi\`{e}re, \emph{M\'{e}thodes semi-classiques et th\'{e}orie
		spectrale},
	https://www-fourier.univ-grenoble-alpes.fr/~ycolver/All-Articles/93b.pdf.
	
	\bibitem[CDV79]{CdVhelv}
	Y.~Colin De~Verdi\`ere, \emph{Sur le spectre des op\'erateurs elliptiques \`a
		bicaract\'eristiques toutes p\'eriodiques}, Comment. Math. Helv. \textbf{54}
	(1979), no.~3, 508--522. \MR{543346}
	
	\bibitem[CdV80]{CdV}
	Y.~Colin~de Verdi\`{e}re, \emph{{Spectre conjoint d'op{\'e}rateurs
			pseudo-diff{\'e}rentiels qui commutent. II. Le cas intégrable}},
	Mathematische Zeitschrift \textbf{171} (1980), 51--74.
	
	\bibitem[Cha88]{ChaIsa}
	A.-M. Charbonnel, \emph{Comportement semi-classique du spectre conjoint
		d'op\'{e}rateurs pseudodiff\'{e}rentiels qui commutent}, Asymptotic Anal.
	\textbf{1} (1988), no.~3, 227--261. \MR{962310}
	
	\bibitem[Com88]{combescure}
	M.~Combescure, \emph{The quantum stability problem for some class of
		time-dependent {H}amiltonians}, Ann. Physics \textbf{185} (1988), no.~1,
	86--110. \MR{954669}
	
	\bibitem[DI17]{DelIme}
	J.~Delort and R.~Imekraz, \emph{Long-time existence for the semilinear
		{K}lein-{G}ordon equation on a compact boundary-less {R}iemannian manifold},
	Comm. Partial Differential Equations \textbf{42} (2017), no.~3, 388--416.
	\MR{3620892}
	
	\bibitem[DLSV02]{duclos2}
	P.~Duclos, O.~Lev, P.~{S}{t}ov\'{\i}\v{c}ek, and M.~Vittot, \emph{Weakly
		regular {F}loquet {H}amiltonians with pure point spectrum}, Rev. Math. Phys.
	\textbf{14} (2002), no.~6, 531--568. \MR{1915516}
	
	\bibitem[DS96]{Duclo1}
	P.~Duclos and P.~{S}{t}ov\'{\i}\v{c}ek, \emph{Floquet {H}amiltonians with pure
		point spectrum}, Comm. Math. Phys. \textbf{177} (1996), no.~2, 327--347.
	\MR{1384138}
	
	\bibitem[Fas05]{Fasso}
	F.~Fass\`o, \emph{Superintegrable {H}amiltonian systems: geometry and
		perturbations}, Acta Appl. Math. \textbf{87} (2005), no.~1-3, 93--121.
	\MR{2151125}
	
	\bibitem[FFS15]{FFS}
	F. Fass\`o, D. Fontanari, and D.~A. Sadovski\'{\i},
	\emph{An application of {N}ekhoroshev theory to the study of the perturbed
		hydrogen atom}, Math. Phys. Anal. Geom. \textbf{18} (2015), no.~1, Art. 30,
	23. \MR{3430232}
	
	\bibitem[Fis15]{Fischer}
	V.~Fischer, \emph{Intrinsic pseudo-differential calculi on any compact lie
		group}, Journal of Functional Analysis \textbf{268} (2015), no.~11,
	3404--3477.
	
	\bibitem[GCB16]{GCB}
	M.~Guzzo, L.~Chierchia, and G.~Benettin, \emph{The steep {N}ekhoroshev's
		theorem}, Communications in Mathematical Physics \textbf{342} (2016),
	569--601.
	
	\bibitem[Gio03]{gio_pisa}
	A.~Giorgilli, \emph{Notes on exponential stability of {H}amiltonian systems},
	Pubblicazioni della Classe di Scienze, Scuola Normale Superiore, Pisa. Centro
	di Ricerca Matematica "Ennio De Giorgi" (2003).
	
	\bibitem[GIP09]{GImekraz}
	B.~Gr\'{e}bert, R.~Imekraz, and E.~Paturel, \emph{Normal forms for semilinear
		quantum harmonic oscillators}, Comm. Math. Phys. \textbf{291} (2009), no.~3,
	763--798. \MR{2534791}
	
	\bibitem[Gor70]{Gordon}
	W.~B. Gordon, \emph{On the relation between period and energy in periodic
		dynamical systems}, J. Math. Mech. \textbf{19} (1969/70), 111--114.
	\MR{245930}
	
	\bibitem[H\"o85]{hormander}
	L.~H\"ormander, \emph{The analysis of linear partial differential operators
		i-iii.}, Springer Berlin, 1985.
	
	\bibitem[HR82]{HR82}
	B.~Helffer and D.~Robert, \emph{Propri\'et\'es asymptotiques du spectre
		d'op\'erateurs pseudodiff\'erentiels sur {${\bf R}^{n}$}}, Comm. Partial
	Differential Equations \textbf{7} (1982), no.~7, 795--882. \MR{662451}
	
	\bibitem[Nek72]{NekInt}
	N.~N. Nekhoro\v{s}ev, \emph{Action-angle variables, and their generalizations},
	Trudy Moskov. Mat. Ob\v{s}\v{c}. \textbf{26} (1972), 181--198. \MR{365629}
	
	\bibitem[Nek77]{Nek77}
	N.~N. Nekhoroshev, \emph{Exponential estimate of the stability of near
		integrable {H}amiltonian systems}, Russ. Math. Surveys \textbf{32} (1977),
	no.~6, 1--65.
	
	\bibitem[Nek79]{Nek79}
	\bysame, \emph{An exponential estimate of the time of stability of nearly
		integrable {H}amiltonian systems. {II}}, Trudy Sem. Petrovsk. (1979), no.~5,
	5--50. \MR{549621}
	
	\bibitem[Nek05]{Nek_Gordon}
	\bysame, \emph{Types of integrability on a submanifold and generalizations of
		{G}ordon's theorem}, Tr. Mosk. Mat. Obs. \textbf{66} (2005), 184--262.
	\MR{2193433}
	
	\bibitem[Nie06]{Nied06}
	L.~Niederman, \emph{Hamiltonian stability and subanalytic geometry}, Ann. Inst.
	Fourier (Grenoble) \textbf{56} (2006), no.~3, 795--813. \MR{2244230}
	
	\bibitem[Pro07]{procesi_padre}
	C.~Procesi, \emph{Lie groups}, Universitext, Springer, New York, 2007, An
	approach through invariants and representations. \MR{2265844}
	
	\bibitem[PVN14]{PSVG}
	A.~Pelayo and S.~V{\~u}~Ng{o}c, \emph{Semiclassical inverse spectral theory for
		singularities of focus--focus type}, Commun. Math. Phys. \textbf{329} (2014),
	no.~2, 809--820 (en).
	
	\bibitem[R{\"{u}}s01]{Rus01}
	H.~R{\"{u}}ssmann, \emph{Invariant tori in non-degenerate nearly integrable
		{H}amiltonian systems}, Regul. Chaotic Dyn. \textbf{6} (2001), no.~2,
	119--204. \MR{1843664}
	
\end{thebibliography}

\providecommand{\bysame}{\leavevmode\hbox to3em{\hrulefill}\thinspace}
\providecommand{\MR}{\relax\ifhmode\unskip\space\fi MR }
\providecommand{\MRhref}[2]{%
	\href{http://www.ams.org/mathscinet-getitem?mr=#1}{#2}
}
\providecommand{\href}[2]{#2}

\end{document}